\newcommand{\dd}[1]{\mathrm{d}{#1}}

\documentclass{amsart}

\usepackage{amsmath}
\usepackage{amssymb}
\usepackage{amsfonts}
\usepackage{MnSymbol}
\usepackage{graphics} %% add this and next lines if pictures should be in esp format
\usepackage{epsfig} %For pictures: screened artwork should be set up with an 85 or 100 line screen
\usepackage{psfrag}
\usepackage{pdfsync}
\usepackage[usenames]{color}
\usepackage{cancel}

\definecolor{arancio}{rgb}{0.90,0.50,0.20}

\definecolor{blu}{rgb}{0.,0.,1.}

\definecolor{pavone}{rgb}{0.00,0.00,0.63}

\definecolor{malva}{rgb}{0.10,0.50,0.50}

\definecolor{rosso}{rgb}{1.,0.,0.}

\definecolor{geranio}{rgb}{0.90,0.00,0.20}

\definecolor{cerulean}{rgb}
{0.0, 0.48, 0.65}

\newtheorem{theorem}{Theorem}[section]
\newtheorem{corollary}[theorem]{Corollary}

\newtheorem{lem}[theorem]{Lemma}

\theoremstyle{definition}
\newtheorem{definition}{Definition}[section]
\newtheorem{remark}{Remark}[section]

\newcommand{\ep}{\varepsilon}
\newcommand{\N}{\mathbb{N}}
\newcommand{\R}{\mathbb{R}}

\date{\today}

\newcommand{\bcl}{\begin{center}}
\newcommand{\ecl}{\end{center}}
\newcommand{\brl}{\begin{right}}
\newcommand{\erl}{\end{right}}
\newcommand{\ben}{\begin{enumerate}}
\newcommand{\barr}{\begin{array}}
\newcommand{\earr}{\end{array}}
\newcommand{\btab}{\begin{tabular}}
\newcommand{\etab}{\end{tabular}}
\newcommand{\bdoc}{\begin{document}}
\newcommand{\edoc}{\end{document}}
\newcommand{\beqy}{\begin{eqnarray}}

\newcommand{\beq}{\begin{equation}}
\newcommand{\beqi}{\begin{eqnarray*}}
\newcommand{\bitem}{\begin{itemize}}
\newcommand{\brem}{\begin{remark}}
\newcommand{\erem}{\end{remark}}
\newcommand{\eitem}{\end{itemize}}
\newcommand{\nln}{\newline}
\newcommand{\newt}{\newtheorem}

\renewcommand{\a }{\alpha }
\renewcommand{\b }{\beta }
\newcommand{\g }{\gamma}
\newcommand{\G }{\Gamma }
\renewcommand{\d }{\delta }

\newcommand{\D }{\Delta }
\newcommand{\e }{\epsilon }
\newcommand{\z }{\zeta }
\renewcommand{\l }{\lambda }
\renewcommand{\L }{\Lambda }
\newcommand{\m }{\mu }
\newcommand{\n }{\tau }
\renewcommand{\r }{\rho }
\newcommand{\s }{\sigma }
\newcommand{\Sig }{\Sigma }
\renewcommand{\t }{\tau }
\newcommand{\var }{H }
\renewcommand{\o }{\omega }
\renewcommand{\O }{\Omega }
\newcommand{\eps}{\varepsilon}
\newcommand{\essinf}{\text{ess\,inf\,}}

\newcommand{\supp}{\text{\rm supp}\,}
\newcommand{\sgn}{\text{\rm sgn}\,}
\title[A free boundary problem for binary fluids]
{A free boundary problem for binary fluids}

\author[Benzi]{Roberto Benzi}
\address{Dipartimento di Fisica, Universit\`a di Roma ``Tor Vergata'', 
Via della Ricerca Scientifica, 00133 Roma, Italy} 
\email{roberto.benzi@gmail.com}
\author[Bertsch]{Michiel Bertsch}
\address{Dipartimento di Matematica, Universit\`a di Roma ``Tor Vergata'', 
Via della Ricerca Scientifica, 00133 Roma, Italy \\ and
Istituto per le Applicazioni del Calcolo "M. Picone", CNR, Roma, Italy} 
\email{bertsch.michiel@gmail.com}
\author[Deangelis]{Francesco Deangelis}
\address{Dipartimento di Matematica, Universit\`a di Roma ``Tor Vergata'', 
Via della Ricerca Scientifica, 00133 Roma, Italy}
\thanks{}
\email{fdeangelis@hotmail.it}



\date{\today}

\begin{document}

\bibliographystyle{h-elsevier2}

\begin{abstract}
A free boundary problem for the dynamics of a glasslike binary fluid naturally leads to a singular perturbation problem 
for a strongly degenerate parabolic partial differential equation in 1D. 
We present a conjecture for an asymptotic formula for the velocity of the free boundary and prove a weak version of the conjecture.
The results are based on the analysis of a family of local travelling wave solutions.
\end{abstract}

\maketitle

\section{Introduction}
In \cite{benzi}, Benzi, Sbragaglia, Bernaschi and Succi propose a phase-field model for binary fluid mixtures.
The fluid mixture is described by an order parameter $\phi$
which, in 1D, satisfies the partial differential equation (PDE) 
\begin{equation}\label{phi1}
\begin{aligned}
\phi_t
&=((D_0+D_2\phi^2)\phi_x)_x -D_2\phi\phi_x^2  +\phi(1-\phi^2)\\
&=(D_0+D_2 \phi^2)\phi_{xx}+D_2\phi\phi_x^2  +\phi(1-\phi^2) \qquad\text{in }Q=(a,b)\times (0,T],
\end{aligned}
\end{equation}
where $D_0$ and $D_2$ are positive constants and $a < b$. The parabolic PDE \eqref{phi1} corresponds to a gradient flow in $L^2$ for the energy functional 
\begin{equation*}
F[\phi]\!=\!\int_a^b\! \left(V(\phi)+\tfrac12 D(\phi) (\phi')^2\right)\dd{x},\qquad V(\phi)\!=\!-\tfrac 12 \phi^2 + \tfrac 14 \phi^2, \quad D(\phi)\!=\!D_0+D_2\phi^2.
\end{equation*}
The double-well potential $V(\phi)$ favours  the two bulk phases, $\phi = \pm 1$.
The stiffness coefficient $D(\phi)$ controls the cost of building and maintaining the interface, $\phi=0$, between the two fluids.

The model is characterised by the choice of $D(\phi)$. 
For small values of $\eps=D_0/D_2$, the mobility of the fluid mixture is relatively small at points where $\phi\approx 0$. 
Physically this leads to a soft-glassy behaviour  of the fluid mixture. In particular, computational evidence in \cite{benzi} suggests 
that the two-fluid interface becomes almost immobile for small values of $\eps$.   
The main purpose of the present paper is to prove this phenomenon analytically and possibly quantify it in terms of the smallness of $\ep$.

Rescaling $x$ by a factor $\sqrt {2/D_2}$ (and changing the interval $(a,b)$ accordingly)
and $t$ by a factor 2, equation $\eqref{phi1}$ becomes
\begin{equation}
\label{eq: evolution equation for phi}
\phi_t =(\varepsilon+ \phi^2)\phi_{xx}+ \phi\phi_x^2  +\tfrac 12 \phi(1-\phi^2) \quad\text{in }Q.
\end{equation}
To understand its mathematical structure for small values of $\eps$  it is convenient to replace $\phi$ by
\begin{equation}
\label{eq: change of variable}
u(x,t) = U_\eps(\phi(x,t)), \quad U_\eps(\phi) = 2\int_0^{\phi} \sqrt{\varepsilon+s^2} \, \dd{s}.
\end{equation}
Since the map $U_\eps$ is strictly increasing, odd and onto in $\R$, its inverse $\Phi_\eps$ is well defined and odd. 
This change of variable transforms \eqref{eq: evolution equation for phi} into
\begin{equation}
\label{eq: eqn with eps}
u_t = (\varepsilon+\Phi_\varepsilon^2(u))u_{xx} +  \Phi_\varepsilon(u)(1-\Phi_\varepsilon^2(u))\sqrt{\varepsilon+\Phi_\varepsilon^2(u)}\quad\text{in }Q.
\end{equation}
In the formal limit  $\eps \to 0 $, $U_\ep(\phi) \to |\phi|\phi$ and $\Phi_\varepsilon(u) \to u/\sqrt{|u|}$. The limit equation is 
\begin{equation}
\label{eq: eqn limit pbm}
u_t = |u|u_{xx}+u(1-|u|)\quad\text{in }Q.
\end{equation}

The parabolicity of equation \eqref{eq: eqn limit pbm} is strongly degenerate. A few decades ago it was studied, 
mostly without reaction term, in the context of {\it nonnegative} solutions 
\cite{bertsch-discontinous, bertsch-nonuniqueness, bertsch-positivity, dal-passo-degenerate, ughi}. 
Several singular phenomena were identified which indicate that its degeneracy is indeed much stronger 
than that of, for example, the well-known porous medium 
and $p-$laplacian equations \cite{DiBenedetto, vazquez-pme}. 
We mention a few of them: solutions are not always uniquely determined by their initial-boundary data and the support of nonnegative solutions is 
non-expanding in time and may also shrink. In higher spatial dimension solutions of a slightly different but similar equation may even be 
discontinuous.

As far as we are aware of, sign-changing solutions of \eqref{eq: eqn limit pbm} were never studied analytically. 
In this paper we consider a class of classical solutions $u_\ep$ of the uniformly parabolic equation \eqref{eq: eqn with eps} 
which initially have a finite number of interfaces,
and we show that for vanishing $\ep$ they 
converge  to a well-defined solution $u$ of the degenerate parabolic limit equation \eqref{eq: eqn limit pbm} and
the interfaces of $u$ are constant in time. In particular, we prove the phenomenon which was numerically observed in \cite{benzi},
namely that for positive but small values of $\ep$ the interfaces of $u_\ep$  are almost immobile.

To quantify the latter result we analyse the existence of a family of local {\it travelling wave solutions} of \eqref{eq: eqn with eps}.
As we shall see, this naturally leads to the following conjecture for an asymptotic formula of the velocity 
of the interface
$x=\zeta_\eps(t)$
(here for simplicity we consider the case that $u_\eps(x,t)$ is monotonic in $x$, so the interface is unique):
\begin{equation}\label{conjecture}
\zeta_\ep'(t)=\frac{u_x(x_1^+,t)-u_x(x_1^-,t)}{2\log \eps}(1+o(1))\quad\text{as }\eps\to 0.
\end{equation}
Here $x_1$ is the position of the constant interface 
of the limit solution $u$. The travelling wave solutions do satisfy \eqref{conjecture}, and in that case the wave velocity  and the right-hand side 
are independent of time.
We also prove that the one-sided spatial derivatives of $u$ at $x_1$ exist and are continuous with respect to $t$ (with the 
possible exception of at most two values of $t$, see Theorem \ref{main result 2});
generically they do not not coincide, as was already observed in \cite{benzi} in the case of stationary solutions
(in \cite{benzi} they were referred to as {\it compactons}).
For the moment being we are not able to prove the conjecture, but we do prove a weaker version of \eqref{conjecture} (see Theorem \ref{main result 3}
and the discussion in Section \ref{proof of main results}).

The paper is organised as follows. 
In Section \ref{limit problem} we collect some preliminary results on the limit problem which will be proved in Section \ref{pf limit pb}. 
In Section \ref{main results} we present the main results and in 
Section \ref{TWs} we analyse the family of travelling wave solutions.
In Section \ref{convergence to the limit problem} we characterise the limit of solutions of \eqref{eq: eqn with eps} for vanishing $\eps$
and show that asymptotically, for vanishing $\eps$, interfaces do not move.
In Section \ref{section regularity} we prove the regularity result for the one-sided spatial derivatives of the limit solution, and
in Section \ref{proof of main results} we prove the weak version of 
the interface condition \eqref{conjecture}.

%%%%%%%%%%%%%%%%%%%%%%%%%%%%%%%%%%%%%%%%%%%%%%%%%%%%%
%%%%%%%%%%%%%%%%%%%%%%%%%%%%%%%%%%%%%%%%%%%%%%%%%%%%%
%%%%%%%%%%%%%%%%%%%%%%%%%%%%%%%%%%%%%%%%%%%%%%%%%%%%%
%%%%%%%%%%%%%%%%%%%%%%%%%%%%%%%%%%%%%%%%%%%%%%%%%%%%%
%%%%%%%%%%%%%%%%%%%%%%%%%%%%%%%%%%%%%%%%%%%%%%%%%%%%%
%%%%%%%%%%%%%%%%%%%%%%%%%%%%%%%%%%%%%%%%%%%%%%%%%%%%%
%%%%%%%%%%%%%%%%%%%%%%%%%%%%%%%%%%%%%%%%%%%%%%%%%%%%%
%%%%%%%%%%%%%%%%%%%%%%%%%%%%%%%%%%%%%%%%%%%%%%%%%%%%%
%%%%%%%%%%%%%%%%%%%%%%%%%%%%%%%%%%%%%%%%%%%%%%%%%%%%%
%%%%%%%%%%%%%%%%%%%%%%%%%%%%%%%%%%%%%%%%%%%%%%%%%%%%%
%%%%%%%%%%%%%%%%%%%%%%%%%%%%%%%%%%%%%%%%%%%%%%%%%%%%%
%%%%%%%%%%%%%%%%%%%%%%%%%%%%%%%%%%%%%%%%%%%%%%%%%%%%%
%%%%%%%%%%%%%%%%%%%%%%%%%%%%%%%%%%%%%%%%%%%%%%%%%%%%%

\section{The limit problem}\label{limit problem}

We consider the problem for equation \eqref{eq: eqn limit pbm} in a bounded interval $(a,b)$, and impose that $u$ is in one of phases $u=\pm 1$ 
at $a$ and $b$:
\begin{equation}\label{limit Dirichlet problem}
\begin{cases}
%u_t = (\varepsilon+\Phi_\varepsilon^2(u))u_{xx} +  \Phi_\varepsilon(u)(1-\Phi_\varepsilon^2(u))\sqrt{\varepsilon+\Phi_\varepsilon^2(u)}&\text{in }Q=(a,b)\times (0,T]\\
u_t = |u| u_{xx} + u(1-|u|) &\text{in }Q=(a,b)\times (0,T]\\
u(a,t)=-1,\ u(b,t)=1 &\text{for }t\in(0,T]\\
u(x,0)=u_0(x) &\text{for }x\in (a,b).
\end{cases}
\end{equation}
Here $u_0:[a,b]\to \R$ is a given initial function which satisfies 
\begin{equation}\label{hyp u_0}
\begin{cases}
u_0\in C([a,b]);
\ u_0(a)=-1,\, u_0(b)=1;\ \\
\text{the number of zeros of $u_0$ in $(a,b)$ is finite};\\ 
\text{$u_0$ changes sign at each of its zeros}. 
\end{cases}
\end{equation}

\begin{definition}\label{weak solution}
Let $u \in L^\infty(Q)\cap L^2(0,T; H^1(a,b))$; $u$ is called a {\it weak solution} of problem \eqref{limit Dirichlet problem} if  $u(a,t)=-1$ and $u(b,t)=1$ 
for a.e.~$t\in (0,T)$, and
\begin{equation*}
\int_a^bu_0(x) \psi(x,0)\dd x + \iint_{Q}\left(u \psi_t- |u|u_x  \psi_x  -  u_x^2 \,\sgn \!(u)\psi +u(1-|u|)\psi\right) \dd x\dd t = 0 
\end{equation*} 
for all $\psi \in C^{1,1}_c((a,b)\times [0,T))$.
\end{definition}

The strong degeneracy of the parabolicity causes singular phenomena which, since a few decades, are known  for nonnegative solutions.
In particular, nonnegative weak solutions have non-expanding spatial supports and suffer various nonuniqueness phenomena 
(\cite{bertsch-nonuniqueness, dal-passo-degenerate, ughi}). There do exist nonnegative weak solutions with spatial supports which 
are independent of time. 

To understand the case of sign-changing solutions we consider initial data which satisfy \eqref{hyp u_0}. 
To be more specific, let $x_1<x_2<\dots<x_k$ be the finite 
number of zeros of $u_0$ in $(a,b)$: 
\begin{equation}\label{zero's}
\mathcal N_0:=\{x_1,x_2,\dots,x_k\}, \qquad \begin{cases} u_0>0 \text{ in $(x_i,x_{i+1})$ if $i$ is odd}\\u_0<0 \text{ in $(x_i,x_{i+1})$ if $i$ is even}\end{cases}
\end{equation}
(obviously $u_0<0$ in $[a,x_1)$ and $u_0>0$ in $(x_k,b]$).
The idea is that we can solve the problem for $u$ independently in each interval $(x_i,x_{i+1})$ with homogeneous 
Dirichlet data at $x_i$ and $x_{i+1}$, without creating zeros of $u$ at the interior of the interval. This naturally leads to the concept of  ``classical'' 
solution of the equation.

\begin{definition}\label{classical solution} Let \eqref{hyp u_0} and \eqref{zero's} be satisfied, and set 
\begin{equation}\label{mathcal N}
\mathcal N:=\mathcal N_0\times [0,T].
\end{equation}
A function $u\in C(\overline Q)$ is called a {\it classical solution} of problem \eqref{limit Dirichlet problem} if   
\begin{itemize}
\item[-] $u(x,t)=0$ if and only if  $(x,t)\in \mathcal N$; 
\item[-] $u\in C^{2,1}(Q\setminus \mathcal N)$; 
\item[-] the second and third equation of  \eqref{limit Dirichlet problem} are satisfied; 
\item[-] the first equation of \eqref{limit Dirichlet problem} is satisfied in $Q\setminus \mathcal N$.
\end{itemize}
\end{definition}

Observe that the continuity of $u$ and hypothesis \eqref{hyp u_0} on $u_0$ imply that, for all $t\in (0, T]$, $u(\cdot,t)\ne 0$ in $( x_i,x_{i+1})$, and 
$u(\cdot,t)$ has the same sign as $u_0$ in $(x_i,x_{i+1})$.
The same observation applies to $(a,x_1)$ and $(x_k,b)$. 

Problem \eqref{limit Dirichlet problem} is well-posed in the class of classical solutions, and the class of classical solutions is a uniqueness class in the set of weak solutions:

\begin{theorem}
\label{well-posedness limit problem}
Let $u_0$ satisfy \eqref{hyp u_0} and \eqref{zero's}.
Then problem \eqref{limit Dirichlet problem} admits a unique classical solution $u \in C^{2,1}(Q) \cap C(\overline{Q})$. 
In addition $u$ is also a weak solution of problem \eqref{limit Dirichlet problem}.
\end{theorem}
\medskip

The proof of Theorem~\ref{well-posedness limit problem} is based on standard techniques and in Section \ref{pf limit pb} we sketch its main lines.
Below (Theorem \ref{main result 2} and Section \ref{section regularity}) we shall establish additional regularity properties of $u$.

%%%%%%%%%%%%%%%%%%%%%%%%%%%%%%%%%%%%%%%%%%%%%%%%%%%%%
%%%%%%%%%%%%%%%%%%%%%%%%%%%%%%%%%%%%%%%%%%%%%%%%%%%%%
%%%%%%%%%%%%%%%%%%%%%%%%%%%%%%%%%%%%%%%%%%%%%%%%%%%%%
%%%%%%%%%%%%%%%%%%%%%%%%%%%%%%%%%%%%%%%%%%%%%%%%%%%%%
%%%%%%%%%%%%%%%%%%%%%%%%%%%%%%%%%%%%%%%%%%%%%%%%%%%%%
%%%%%%%%%%%%%%%%%%%%%%%%%%%%%%%%%%%%%%%%%%%%%%%%%%%%%
%%%%%%%%%%%%%%%%%%%%%%%%%%%%%%%%%%%%%%%%%%%%%%%%%%%%%
%%%%%%%%%%%%%%%%%%%%%%%%%%%%%%%%%%%%%%%%%%%%%%%%%%%%%
%%%%%%%%%%%%%%%%%%%%%%%%%%%%%%%%%%%%%%%%%%%%%%%%%%%%%
%%%%%%%%%%%%%%%%%%%%%%%%%%%%%%%%%%%%%%%%%%%%%%%%%%%%%
%%%%%%%%%%%%%%%%%%%%%%%%%%%%%%%%%%%%%%%%%%%%%%%%%%%%%
%%%%%%%%%%%%%%%%%%%%%%%%%%%%%%%%%%%%%%%%%%%%%%%%%%%%%
%%%%%%%%%%%%%%%%%%%%%%%%%%%%%%%%%%%%%%%%%%%%%%%%%%%%%

\section{Main results}\label{main results}
The main results of the paper concern the behaviour of the unique solution 
$u_\eps\in C^{2,1}(Q)\cap C(\overline Q)$
of the problem
\begin{equation}\label{problem with eps}
\begin{cases}
u_t = (\varepsilon+\Phi_\varepsilon^2(u))u_{xx} +  \Phi_\varepsilon(u)(1-\Phi_\varepsilon^2(u))\sqrt{\varepsilon+\Phi_\varepsilon^2(u)}&\text{in }Q\\
u(a,t)=-u_{1\eps},\ u(b,t)=u_{1\eps} &\text{for }t\in(0,T]\\
u(x,0)=u_{0\eps}(x) &\text{for }x\in (a,b)
\end{cases}
\end{equation}
for small values of $\eps>0$. Here $\Phi_\ep$ is defined as in the Introduction, $a<b$, $T>0$ and $u_{1\eps}>0$ are constants such that 
\begin{equation}\label{boundary condition eps}
\Phi_\eps (u_{1\eps})=1 \qquad (\Rightarrow \Phi_\eps ( -u_{1\eps})=-1, \text{ and } u_{1\eps}\to 1 \text{ as }\eps\to 0),
\end{equation}
and the initial functions $u_{0\eps}$ satisfy
\begin{equation}\label{u_0 eps}
\begin{cases}
u_{0\eps}\in C^\infty ([a,b]), \quad  u_{0\eps}(a)=-u_{1\eps}, \quad  u_{0\eps}(b)=u_{1\eps}\\
u_{0\eps}(x)=0 \Leftrightarrow x\in \mathcal N_0, \quad u_{0\eps}\to u_0\text{ in $C([a,b])$ as }\eps \to 0.
\end{cases}
\end{equation}

The first result shows that $u_\eps$ converges uniformly to the solution of the limit problem for vanishing $\eps$, which implies that away from the set $\mathcal N$,
defined by \eqref{mathcal N}, $u_\eps(x,t)$ has the same sign as $u_0(x)$ if $\eps$ is small enough. More precisely we have:

\begin{theorem}\label{main result 1} Let $u_0$ satisfy \eqref{hyp u_0} and \eqref{zero's},  let $\eps>0$ and let $u_{0\eps}$ satisfy \eqref{u_0 eps}.
Let $u_\eps\in C^{2,1}(Q)\cap C(\overline Q)$ be the solution of problem \eqref{problem with eps}, $u$ the unique
classical solution of the limit problem 
\eqref{limit Dirichlet problem} defined by Theorem \ref{well-posedness limit problem}, 
and $\mathcal N$ the set defined by \eqref{mathcal N}. 
Then $u_\eps \to u$ uniformly in $Q$ and in $C^{2,1}_{\rm loc}(Q\setminus \mathcal N)$, and 
\begin{equation}\label{interfaces do not move}
\sup_{\{(x,t)\in Q;\, u_\eps(x,t)= 0\}}\text{\rm distance\,}( (x,t),\mathcal N)\to 0 \quad\text{as }\eps \to 0.
\end{equation}
\end{theorem}

To state the asymptotic result on the free boundary condition  \eqref{conjecture} 
we need the following regularity result for the one-sided spatial derivatives of the limit problem.

\begin{theorem}\label{main result 2} Let $u_0$ satisfy \eqref{hyp u_0} and let $u$ be the unique
classical solution of the limit problem \eqref{limit Dirichlet problem}, defined by Theorem \ref{well-posedness limit problem}. 
Then 
\begin{itemize}
\item[$(i)$]  $u_x\in L^\infty(((a,b)\setminus \mathcal N_0)\times (t_0,T))$ for all $t_0\in (0,T)$,
and, for all $t\in (0,T]$, the function $x\mapsto u_x(x,t)$ has at most a jump discontinuity at $x_i\in \mathcal N_0$ $(i=1,\dots,k)$;
\item[$(ii)$] for all $i=1,\dots,k$ there exist $\tau_i^\pm\in [0,\infty]$ $(\text{independent of }T)$ such that 
$$
u_x(x_i^\pm,t)\begin{cases} =0 &\text{if } 0<t<\tau_i^\pm \\ \ne 0&\text{if }t>\tau_i^\pm;\end{cases}
$$
in addition the functions $t\mapsto u_x(x_i^\pm,t)$ are continuous in  $(0,T]\setminus \{\tau_i^\pm\}$.
\end{itemize}
\end{theorem}

Concerning  the asymptotic expansion  \eqref{conjecture} we limit ourselves,
for the sake of simplicity, to the case of a strictly increasing initial function:
\begin{equation}\label{hyp u_0 increasing}
u_0\in C([a,b]), \quad u_0 \text{ is strictly increasing in }[a,b],\quad u_0(a)=-1,\, u_0(b)=1.
\end{equation}   
The condition on the approximating initial data $u_{0\eps}$ ($\eps>0$) is changed accordingly:
\begin{equation}\label{hyp u_0,eps increasing}
\begin{cases} 
u_{0\eps}\in C^\infty([a,b]), \quad u'_{0\eps}>0 \text{ in }[a,b],\quad u_{0\eps}(a)=-u_{1\eps},\, u_{0\eps}(b)=u_{1\eps}\\
u_0(x_1)=0\Rightarrow u_{0\eps}(x_1)=0, \quad u_{0\eps}\to u_0 \text{ in } C([a,b])\text{ as } \eps\to 0.
\end{cases}
\end{equation}

\begin{theorem}\label{main result 3} Let $u_0$ satisfy \eqref{hyp u_0 increasing}, let $\eps>0$ and let $u_{0\eps}$ satisfy \eqref{hyp u_0,eps increasing}.
Let $u_\eps\in C^{2,1}(Q)\cap C(\overline Q)$ be the solution of problem \eqref{problem with eps} and let $u$ be the unique
classical solution of the limit problem \eqref{limit Dirichlet problem}, defined by Theorem \ref{well-posedness limit problem}.
Then $u_{\eps x}>0$ in $Q$. Let $x=X_\eps(u,t)$ be defined by $u_\eps(X_\eps(u,t),t)=u$.
Then there exists $0<\delta_\eps\to 0$ as $\eps\to 0$ such that  for all $t\in [0,T]$
\begin{equation}\label{FBC averaged}
\left( \int_{-\delta_\eps}^{\delta_\eps} \frac{\dd{u}}{\varepsilon+\Phi_\varepsilon^2(u)} \right)^{-1} \!\!\! \int_{-\delta_\eps}^{\delta_\eps}  \frac{ X_{\eps t}(u,t)}{\varepsilon+\Phi_\varepsilon^2(u)} \dd{u}
 =  \frac{u_x(x_1^+,t)-u_x(x_1^-,t)}{2\log \eps}(1+o(1)) \quad\text{as }\eps\to0.
\end{equation}
\end{theorem}

The interpretation is immediate: let $x=\zeta_\eps(t)$ be the interface of $u_\eps(x,t)$; then $\zeta_\eps'(t)=X_{\eps t}(0,t)$, so \eqref{FBC averaged} is nothing else 
than the interface
condition \eqref{conjecture} with the left-hand side, $\zeta_\eps'(t)$,  replaced by a weighted average of $X_{\eps t}(u,t)$ in a neighbourhood $(-\delta_\eps,\delta_\eps)$ 
of $u=0$, a neighbourhood which shrinks to a single point as $\eps\to 0$.
    
In Section \ref{proof of main results} (see Remark \ref{possible weak formulations}) we shall briefly discuss a different weak version of \eqref{conjecture}.

%%%%%%%%%%%%%%%%%%%%%%%%%%%%%%%%%%%%%%%%%%%%%%%%%%%%%
%%%%%%%%%%%%%%%%%%%%%%%%%%%%%%%%%%%%%%%%%%%%%%%%%%%%%
%%%%%%%%%%%%%%%%%%%%%%%%%%%%%%%%%%%%%%%%%%%%%%%%%%%%%
%%%%%%%%%%%%%%%%%%%%%%%%%%%%%%%%%%%%%%%%%%%%%%%%%%%%%
%%%%%%%%%%%%%%%%%%%%%%%%%%%%%%%%%%%%%%%%%%%%%%%%%%%%%
%%%%%%%%%%%%%%%%%%%%%%%%%%%%%%%%%%%%%%%%%%%%%%%%%%%%%
%%%%%%%%%%%%%%%%%%%%%%%%%%%%%%%%%%%%%%%%%%%%%%%%%%%%%
%%%%%%%%%%%%%%%%%%%%%%%%%%%%%%%%%%%%%%%%%%%%%%%%%%%%%
%%%%%%%%%%%%%%%%%%%%%%%%%%%%%%%%%%%%%%%%%%%%%%%%%%%%%
%%%%%%%%%%%%%%%%%%%%%%%%%%%%%%%%%%%%%%%%%%%%%%%%%%%%%
%%%%%%%%%%%%%%%%%%%%%%%%%%%%%%%%%%%%%%%%%%%%%%%%%%%%%
%%%%%%%%%%%%%%%%%%%%%%%%%%%%%%%%%%%%%%%%%%%%%%%%%%%%%
%%%%%%%%%%%%%%%%%%%%%%%%%%%%%%%%%%%%%%%%%%%%%%%%%%%%%

\section{Travelling waves}\label{TWs}

In this section we analyse a family of travelling wave solutions (TWs) which play a key role in the proof of Theorem \ref{main result 1} and the
formulation of the conjecture \eqref{conjecture}.

Let $u(x,t)=w(x-ct)$ be a travelling wave solution of \eqref{eq: eqn with eps} with velocity  $c\in\R$:
\begin{equation}\label{TW equation}
-cw'=(\varepsilon+ \Phi_\varepsilon^2(w))w''+\Phi_\varepsilon(w)(1-\Phi_\varepsilon^2(w))\sqrt{\varepsilon+\Phi_\varepsilon^2(w)}.
\end{equation}
Assuming monotonicity of $w(z)$ in an interval, we use 
$w$ as independent variable to reduce the order of this autonomous ODE: 
the function  $p(w)$, defined  by 
$$
p(w(z))=w'(z),
$$
satisfies
\begin{equation}\label{p'}
p'=-\frac{c}{\varepsilon+ \Phi_\varepsilon^2(w)}-\frac {\Phi_\varepsilon(w)(1-\Phi_\varepsilon^2(w))}{p\sqrt{\varepsilon+\Phi_\varepsilon^2(w)}}.
\end{equation}

We are interested in local TWs in a neighbourhood of the ``interface'' $w=0$ for which $w'(z)$ is strictly positive. 
If $p$ is strictly positive, the second term on the right hand side of \eqref{p'} is bounded. Therefore we 
focus on the first term on the RHS which becomes singular at $w=0$ for vanishing $\eps$. Given 
$\delta>0$, we have that
$$
I_\delta:=- c\int_0^\delta \frac{1}{\varepsilon+ \Phi_\varepsilon^2(w)}\,dw
=-2c\int_{0}^{\Phi_\varepsilon(\delta)}\frac 1{\sqrt{\varepsilon+\phi^2}}\,d\phi.
$$
We set $\phi=\sqrt\varepsilon \sinh y$. Since $d\phi =\sqrt\varepsilon \cosh y\,dy$, we obtain that
$$
\begin{aligned}
\int \frac1{\sqrt {\varepsilon+ \phi^2}}\,d\phi
&=\int \frac{\cosh y}{\sqrt {1+ \sinh^2y}}\,dy
= y 
= \log\left(\frac{1}{\sqrt \varepsilon}\left(\phi+\sqrt{\varepsilon+\phi^2}\right)\right),\\
\end{aligned}
$$
and
\begin{equation}\label{integral}
I_\delta=-2c\log\left(\frac{1}{\sqrt \varepsilon}\left(\Phi_\varepsilon(\delta)+\sqrt{\varepsilon+\Phi_\varepsilon^2(\delta)}\right)\right).
\end{equation}
If $\frac{\Phi_\varepsilon(\delta)}{\sqrt{\varepsilon}}$ is large, then
$$
I_\delta\approx-2c\log\left(\frac{2\Phi_\varepsilon(\delta)}{\sqrt{\varepsilon}}\right)
=-2c\left(\log\left(2\Phi_\varepsilon(\delta)\right)-\tfrac 12 \log \varepsilon \right)\approx c\log \varepsilon.
$$

So if $I_\delta$  is bounded away from zero, $c$ vanishes as $\varepsilon \to 0$, and 
since $I_\delta$ represents a variation in $p$, this simple calculation suggests the following result.
Before stating it we observe that, given $B>0$, the function $w_B\in C([0,\infty))$ defined by 
\begin{equation}\label{w_B}
w_B(x)=\max\{B \sinh x -\cosh x  +1,0\} \quad \text{for } x \geq 0,
\end{equation}
%(here $[y]_+=\max\{y,0\}$ for $y\in \R$)  
is a nonnegative steady state of the limit equation \eqref{eq: eqn limit pbm} in $[0,\infty)$ which satisfies
$w_B(0)=0$, $w_B'(0)=B$ and
$$
w_B(x)>0 \quad \text{if }
\begin{cases}
0<x<\log\frac {1+B}{1-B} &\text{if }0<B<1\\
x>0 &\text{if }B\ge 1.
\end{cases}
$$

\begin{lem}
\label{lmm: lemma travelling waves approximation}
Let $B$ and $B_0$ be positive constants. Let 
$w_B\in C([0,\infty))$ be the steady state defined by \eqref{w_B}.
Let $\eps>0$ and set
\begin{equation}\label{estimate c 1}
c_\varepsilon=\frac{B-B_0}{\log\varepsilon}.
\end{equation}
Let $w_{B,\varepsilon}$ be the local solution of the shooting problem
$$
\begin{cases}
-c_\eps w'=(\varepsilon+ \Phi_\varepsilon^2(w))w''+\Phi_\varepsilon(w)(1-\Phi_\varepsilon^2(w))\sqrt{\varepsilon+\Phi_\varepsilon^2(w)} &\text{for }x>0\\
w(0)=0,\ w'(0)=B_0,
\end{cases}
$$
which can be continued as long as it stays bounded.
Then 
\begin{equation}\label{convergence to w_B}
w_{B,\varepsilon}\to w_B \quad\text{in }
\begin{cases} 
C^1_{\rm loc}([0,\log\frac {1+B}{1-B})) &\text{if } 0<B<1\\ 
C^1_{\rm loc}([0,\infty)) &\text{if } B>1.
\end{cases}
\end{equation}
\end{lem}

\begin{proof}
Since $B_0>0$, $w'_{B,\varepsilon}>0$ near $x=0$. As long as $w_{B,\varepsilon}$ remains increasing and bounded, we argue as above and introduce $p_\eps(w)$,
which locally is a solution of
\begin{equation}\label{pb p 1}
\begin{cases}
p'=-\dfrac{c_\varepsilon}{\varepsilon+ \Phi_\varepsilon^2(w)}-\dfrac {\Phi_\varepsilon(w)(1-\Phi_\varepsilon^2(w))}{p\sqrt{\varepsilon+\Phi_\varepsilon^2(w)}} 
&\text{for }w>0 \\
p(0)=B_0.
\end{cases}
\end{equation}
To understand the behaviour of $p_\eps$ near $w=0$ we change variable and set
$$
\begin{aligned}
q_\eps(w)&=p_\eps(w)+\int_{0}^w\dfrac{c_\varepsilon}{\varepsilon+ \Phi_\varepsilon^2(s)}\,ds \\
&=p_\eps(w)+\frac{2(B-B_0)}{\log\varepsilon}\log\left[\tfrac{\Phi_\varepsilon(w)+\sqrt{\varepsilon+\Phi^2_\varepsilon(w)}}
{\sqrt{\varepsilon}}\right],
\end{aligned}
$$
where we have used \eqref{integral}.  The equation for $q_\eps$ is
$$
q_\eps'(w)
=-\dfrac {\Phi_\varepsilon(w)(1-\Phi_\varepsilon^2(w))}{\left(q_\eps(w)- 
\frac{2(B-B_0)}{\log\varepsilon}\log\left[\frac{\Phi_\varepsilon(w)+\sqrt{\varepsilon+\Phi^2_\varepsilon(w)}}
{\sqrt{\varepsilon}}\right]\right)\sqrt{\varepsilon+\Phi_\varepsilon^2(w)}}. 
$$
Let $w_0>0$ be a small number to be chosen below. If $0\le w\le w_0$, then
$$
\begin{aligned}
1\le&\frac{\Phi_\varepsilon(w)+\sqrt{\varepsilon+\Phi^2_\varepsilon(w)}}
{\sqrt{\varepsilon}} 
\le\frac{\Phi_\varepsilon(w_0)+\sqrt{\varepsilon+\Phi^2_\varepsilon(w_0)}}
{\sqrt{\varepsilon}}=
\frac{2 \sqrt{w_0}}{\sqrt{\eps}}(1+o(1)) \\
\end{aligned}
$$
as $\eps\to 0$, whence 
$$
-1+o(1) \le \frac 2{\log\varepsilon}\log\left[\frac{\Phi_\varepsilon(w)+\sqrt{\varepsilon+\Phi^2_\varepsilon(w)}}
{\sqrt{\varepsilon}}\right]\le 0 \quad\text{as }\eps\to 0.
$$ 
Since the factor $\tfrac{\Phi_\eps(w)(1-\Phi_\eps^2(w))}{ \sqrt{\eps + \Phi_\eps^2(w)}}$ 
in the equation for $q_\eps$ is uniformly bounded,
this implies that $q_\eps$ can be made arbitrarily near to $B_0>0$ in the interval $(0,w_0)$ 
by choosing $w_0$ and $\eps$ small enough. 

If $B\ge B_0$ this means that $p_\eps\ge q_\eps$ is bounded away from 0 in $(0,w_0)$ if $w_0$ and $\eps$ 
are chosen small enough. If instead $B<B_0$, it is enough to slightly refine this argument and use that 
$$
(B-B_0)(1+o(1)) \le -\frac{2(B-B_0)}{\log\varepsilon}\log\left[\frac{\Phi_\varepsilon(w)+\sqrt{\varepsilon+\Phi^2_\varepsilon(w)}}
{\sqrt{\varepsilon}}\right]\le 0 \quad\text{as }\eps\to 0
$$
uniformly in $(0,w_0)$, which implies that we can choose $w_0$ and $\eps$ so small that in $(0,w_0)$ the solution $p_\eps$ is contained in an interval which is 
only slightly larger than $[B,B_0]$ (and in particular bounded away from 0).

We check that $p_\eps$ changes by the ``right" amount in the interval $(0, w_0)$, namely that $p_\eps(w_0)$ can be made
arbitrarily  close to $B$ by choosing $w_0$ and $\eps$ small enough. Since $q_\eps$ can be made arbitrarily close to $B_0$ in   $(0,w_0)$, 
it is enough to observe that, as $\eps\to 0$, 
$$
-\int_{0}^{w_0}\dfrac{c_\varepsilon}{\varepsilon+ \Phi_\varepsilon^2(w)}\,dw 
= -\frac{2(B-B_0)}{\log\varepsilon}\log\left[\frac{\Phi_\varepsilon(w_0)+\sqrt{\varepsilon+\Phi^2_\varepsilon(w_0)}}
{\sqrt{\varepsilon}}\right]
\to B-B_0. 
$$

Turning to the original variables $w_{B,\eps}(x)$, it follows that the slope $w'_{B,\eps}(x)$ changes in a small right neighbourhood of the origin from $B_0$ 
to approximately $B$.
For larger values of $x$ we then use that $w_{B,\eps}$ depends continuously on $\eps$ as long as it stays bounded and positive. This easily leads
to \eqref{convergence to w_B}.
\end{proof}

The above result has its natural counterpart for {\it nonpositive} steady states $\tilde w_A$ of the limit equation with support in $(-\infty,0]$,
\begin{equation}\label{w_A}
\tilde w_A(x)=\min\{A \sinh x+\cosh x -1,0\}  \quad \text{for } x \le 0.
\end{equation}

\begin{corollary}
\label{lmm: lemma travelling waves approximation bis}
Let $A$ and $A_0$ be positive constants. Let 
$\tilde w_A\in C((-\infty,0])$ be the steady state defined by \eqref{w_A}.
Let $\eps>0$ and set
\begin{equation}\label{estimate c 1 bis}
\tilde c_\varepsilon=-\frac{A-A_0}{\log\varepsilon}.
\end{equation}
Let $\tilde w_{A,\varepsilon}$ be the local solution of the shooting problem
$$
\begin{cases}
-\tilde c_\eps w'=(\varepsilon+ \Phi_\varepsilon^2(w))w''+\Phi_\varepsilon(w)(1-\Phi_\varepsilon^2(w))\sqrt{\varepsilon+\Phi_\varepsilon^2(w)} &\text{for }x<0\\
w(0)=0,\ w'(0)=A_0,
\end{cases}
$$
which can be continued as long as it stays bounded.
Then 
\begin{equation}\label{convergence to w_A}
\tilde w_{A,\varepsilon}\to \tilde w_A \quad\text{in }
\begin{cases} 
C_{\rm loc}((-\log\frac {1+A}{1-A},0]) &\text{if } 0<A<1\\ 
C_{\rm loc}((-\infty,0]) &\text{if } A>1.
\end{cases}
\end{equation}
\end{corollary}

We can easily merge the approximations of a nonnegative steady state in $[0,\infty)$ and a nonpositive one in $(-\infty,0]$ 
to construct an approximation $w_{AB,\eps}$ of the following steady state with changing sign: 
\begin{equation}\label{w_AB}
w_{AB}(x)=
\begin{cases}
\tilde w_A(x)= \min\{A \sinh x+\cosh x -1,0\}&\text{if }x\le 0\\
w_B(x)= \max\{B \sinh x -\cosh x  +1,0\} &\text{if }x> 0.\\
\end{cases}
\end{equation}
To do so we require that the two wave velocities coincide, $c_\eps=\tilde c_\eps$, as well as the two shooting parameters, $A_0=B_0$. This means that
$$
A_0=B_0=\tfrac 12 (A+B), \qquad c_\eps= \tilde c_\eps =\frac{B-A}{2\log\varepsilon}.
$$
Combining Lemma \ref{lmm: lemma travelling waves approximation} and Corollary \ref{lmm: lemma travelling waves approximation bis}, we obtain the following result.

\begin{theorem}
\label{thm: travelling waves approximation}
Let $A$ and $B$ be positive constants and let $\eps\in (0,1]$. Let 
$w_{AB}\in C(\R)$ be the steady state of the limit equation \eqref{eq: eqn limit pbm} defined by \eqref{w_AB}. 
Then there exists a travelling wave solution 
$w_{AB,\varepsilon}(x-c_\varepsilon t)$ of equation \eqref{eq: eqn with eps} with velocity  
\begin{equation}\label{estimate c}
c_\varepsilon=\frac{B-A}{2\log\varepsilon}
\end{equation}
such that $w_{AB,\varepsilon}(0)=0$ and $w_{AB,\varepsilon}\to w_{AB}$ in $C_{\rm loc}(J_{AB})$ as $\varepsilon\to 0$, where 
$$
J_{AB}=
\begin{cases}
\R &\text{if }A,B\ge 1\\
(-\log\frac{1+A}{1-A}, \log\frac{1+B}{1-B})&\text{if }A,B\in (0,1)\\
(-\infty,\log\frac{1+B}{1-B})&\text{if }A\ge 1,\, 0<B< 1\\
(-\log\frac{1+A}{1-A},\infty)&\text{if }0<A<1,\, B\ge  1.\\
\end{cases}
$$ 
\end{theorem}

For later use we observe that if $A>1$, the steady state $\tilde w_A\in C^1((-\infty,0])$, defined by 
\eqref{w_A}, has an inflection point at $-\tfrac 12 \log\frac {A+1}{A-1}$ and 
\begin{equation}\label{w_A'}
\tilde w_A'\ge \tilde w_A'\left(-\tfrac{1}{2}\log\tfrac{A+1}{A-1}\right) =\sqrt{A^2-1} 
\quad\text{in }(-\infty,0].
\end{equation}

%%%%%%%%%%%%%%%%%%%%%%%%%%%%%%%%%%%%%%%%%%%%%%%%%%%%%
%%%%%%%%%%%%%%%%%%%%%%%%%%%%%%%%%%%%%%%%%%%%%%%%%%%%%
%%%%%%%%%%%%%%%%%%%%%%%%%%%%%%%%%%%%%%%%%%%%%%%%%%%%%
%%%%%%%%%%%%%%%%%%%%%%%%%%%%%%%%%%%%%%%%%%%%%%%%%%%%%
%%%%%%%%%%%%%%%%%%%%%%%%%%%%%%%%%%%%%%%%%%%%%%%%%%%%%
%%%%%%%%%%%%%%%%%%%%%%%%%%%%%%%%%%%%%%%%%%%%%%%%%%%%%
%%%%%%%%%%%%%%%%%%%%%%%%%%%%%%%%%%%%%%%%%%%%%%%%%%%%%
%%%%%%%%%%%%%%%%%%%%%%%%%%%%%%%%%%%%%%%%%%%%%%%%%%%%%
%%%%%%%%%%%%%%%%%%%%%%%%%%%%%%%%%%%%%%%%%%%%%%%%%%%%%
%%%%%%%%%%%%%%%%%%%%%%%%%%%%%%%%%%%%%%%%%%%%%%%%%%%%%
%%%%%%%%%%%%%%%%%%%%%%%%%%%%%%%%%%%%%%%%%%%%%%%%%%%%%
%%%%%%%%%%%%%%%%%%%%%%%%%%%%%%%%%%%%%%%%%%%%%%%%%%%%%
%%%%%%%%%%%%%%%%%%%%%%%%%%%%%%%%%%%%%%%%%%%%%%%%%%%%%

\section{Convergence to the limit problem}\label{convergence to the limit problem}

In this section we prove Theorem \ref{main result 1}. 
Let $u_0$ satisfy \eqref{hyp u_0} and \eqref{zero's}, let $u_{0\eps}$ satisfy \eqref{u_0 eps}, %\eqref{boundary condition eps},
let $u_\ep\in C(\overline Q)\cap C^{2,1}(Q)$ be the unique smooth solution of problem \eqref{problem with eps}
and let $u\in C(\overline Q)\cap C^{2,1}(Q\setminus \mathcal N)$ be
the unique classical solution of   the limit problem \eqref{limit Dirichlet problem}.
We recall that $\mathcal N_0\subset (a,b)$ is the set containing the $k$ zeros $x_i$ of $u_0$ and that $\mathcal N=\mathcal N_0\times [0,T]$.

We choose one of the intervals $(x_i,x_{i+1})$. We assume that 
$$
u_0>0\quad\text{in }(x_i,x_{i+1}).
$$ 
The case that $u_0<0$ in $(x_i,x_{i+1})$ can be treated similarly, and
also the intervals $(a,x_1)$ and $(x_k,b)$ are treated in a similar way.
 
We claim that it is enough to prove the following: 
\begin{equation}\label{basic estimate}
\begin{aligned}
&\quad \text{for all $n\in \N$, $n>2(x_{i+1}-x_i)^{-1}$, there exists $\eps_n\in(0,\eps_{n-1})$ such that }\\
&u_\eps(x_i+\tfrac 1n,t)\ge 0\text{ and } u_\eps(x_{i+1}-\tfrac 1n,t)\ge 0 \text{ for }0<t<T \text{ and for all $0<\eps<\eps_n$}.
\end{aligned}
\end{equation}
Indeed, \eqref{basic estimate} and the Comparison Principle in $(x_i+\tfrac 1n, x_{i+1}-\tfrac 1n)\times (0,T)$ imply that
$$
u_\eps\ge \underline v_{\eps, n} \quad\text{in } (x_i+\tfrac 1n, x_{i+1}-\tfrac 1n)\times (0,T) \text{ if }\eps<\eps_n,
$$
where $\underline v_{\eps, n}$ is the unique smooth solution
of the problem
\begin{equation}\label{problem with eps and n}
\begin{cases}
v_t = (\varepsilon+\Phi_\varepsilon^2(v))v_{xx} +  \Phi_\varepsilon(v)(1-\Phi_\varepsilon^2(v))\sqrt{\varepsilon+\Phi_\varepsilon^2(v)}&\text{in }
(x_i+\tfrac 1n, x_{i+1}-\tfrac 1n)\times (0,T]\\
v(x_i+\tfrac 1n,t)=v(x_{i+1}-\tfrac 1n,t)=0 &\text{for }t\in(0,T]\\
v(x,0)=u_{0\eps}(x) &\text{for }x\in (x_i+\tfrac 1n,x_{i+1}-\tfrac 1n)\\
v>0 &\text{in }(x_i+\tfrac 1n, x_{i+1}-\tfrac 1n)\times (0,T].
\end{cases}
\end{equation}
Since $\underline v_{\eps, n+1}\ge \underline v_{\eps,n}>0$ in $(x_i+\tfrac 1n, x_{i+1}-\tfrac 1n)\times (0,T]$, this implies that there exists
$$
\underline v_\eps(x,t)=\lim_{n\to\infty}\underline v_{\eps,n} (x,t) \quad\text{for }(x,t)\in (x_i,x_{i+1})\times [0,T].
$$

On the other hand we use the properties in the intervals $(x_{i-1},x_i)$ and  $(x_{i+1},x_{i+2})$ which are analogous to \eqref{basic estimate} and
which, put together, imply that
$$
\begin{aligned}
&\quad \text{for all sufficiently large $n\in \N$ there exists $\tilde \eps_n\in(0,\eps_{n-1})$ such that }\\
&u_\eps(x_i-\tfrac 1n,t)<0\text{ and } u_\eps(x_{i+1}+\tfrac 1n,t)<0 \text{ for }0<t<T \text{ and for all $0<\eps<\tilde \eps_n$}.
\end{aligned}
$$
Let $\overline v_{\eps, n}$ be the unique smooth solution of problem \eqref{problem with eps and n} with $x_i+\tfrac 1n$ and $x_{i+1}-\tfrac 1n$ replaced
by, respectively, $x_i-\tfrac 1n$ and $x_{i+1}+\tfrac 1n$. Then
$$
u_\eps\le \overline v_{\eps, n} \quad\text{in } (x_i-\tfrac 1n, x_{i+1}+\tfrac 1n)\times (0,T) \text{ if }\eps<\tilde \eps_n,
$$
and since $\overline v_{\eps, n+1}\le \overline v_{\eps,n}$ in $(x_i-\tfrac 1{n+1}, x_{i+1}+\tfrac 1{n+1})\times (0,T]$ we may define
$$
\overline v_\eps(x,t)=\lim_{n\to\infty}\overline v_{\eps,n} (x,t) \quad\text{for }(x,t)\in [x_i,x_{i+1}]\times [0,T].
$$
Hence $0<\underline v_\eps\le \overline v_\eps$ in $(x_i,x_{i+1})\times [0,T]$ and $\overline v_\eps(x_i,t)=\overline v_\eps(x_{i+1},t)=0$. 
By local Schauder type estimates (\cite{LSU}), in  $(x_i,x_{i+1})\times (0,T]$ we may pass to the limit $n\to\infty$ in the equation for 
$\underline v_{\eps,n}$ and $\overline v_{\eps,n}$, whence  both $\underline v_\eps$ and $ \overline v_\eps$ coincide with the unique
solution of the problem 
$$
\begin{cases}
u_t = (\varepsilon+\Phi_\varepsilon^2(u))u_{xx} +  \Phi_\varepsilon(u)(1-\Phi_\varepsilon^2(u))\sqrt{\varepsilon+\Phi_\varepsilon^2(u)}
&\text{in }(x_i, x_{i+1})\times (0,T]\\
u(x_i,t)=u(x_{i+1},t)=0 &\text{for }t\in(0,T]\\
u(x,0)=u_{0\eps}(x) &\text{for }x\in (x_i,x_{i+1})\\
u>0 &\text{in }(x_i, x_{i+1})\times (0,T].
\end{cases}
$$

Now Theorem \ref{main result 1} would follow from a standard diagonal procedure if we could use again local Schauder type estimates in order 
to pass to the limit $\eps\to 0$ and conclude that $\underline v_\eps=\overline v_\eps $
converges to the unique (by Theorem \ref{well-posedness limit problem}) solution of the problem
$$
\begin{cases}
u_t = uu_{xx} +  u(1-u)&\text{in }(x_i, x_{i+1})\times (0,T]\\
u(x_i,t)=u(x_{i+1},t)=0 &\text{for }t\in(0,T]\\
u(x,0)=u_{0}(x) &\text{for }x\in (x_i,x_{i+1})\\
u>0 &\text{in }(x_i, x_{i+1})\times (0,T].
\end{cases}
$$
To justify this, we only need to prove that locally in $(x_i+\tfrac 1n, x_{i+1}-\tfrac 1n)\times [0,T]$ the solution $\underline v_{\eps,n}$ is uniformly bounded away
from 0, which makes the equation for $\underline v_{\eps,n}$ locally uniformly parabolic. 
But establishing such uniform lower bound is easy: given $x_0\in (x_i,x_{i+1})$ there exists a ``small'' steady state
solution $\tilde u(x)$ of the limit equation which is strictly positive in a ``small'' neighbourhood $\mathcal U$ of $x_0$, 
where ``small'' means that $\tilde u\le \tfrac 12u_{0}$ in $\mathcal U$ and $\supp \tilde u\subset (x_i,x_{i+1})$. Arguing as in Section \ref{TWs} it easily follows that 
$\tilde u$ can be approximated by a family of steady state solutions of the equation with $\eps$ which can be used as subsolutions of the problem 
for $\underline v_{\eps,n}$, and since $x_0$ is arbitrary we obtain the uniform lower bound for $\underline v_{\eps,n}$, and thus  for $\underline v_{\eps}$. 
We leave the details to the reader. 

\medskip

It remains to prove \eqref{basic estimate}.
\medskip

\noindent{\it Proof of \eqref{basic estimate}.}
Since $u_{0\eps}\to u_0$ uniformly in $(a,b)$, for all $n\in \N$ there exist $m_n>0$ and $\eps_{1n}>0$ such that 
$\Phi_\varepsilon^2(y)<1$ if $y\in (0,m_n]$ and 
$$
u_{0\eps}>m_n 
\quad\text{in }(x_i+\tfrac 1{2n}, x_{i+1}-\tfrac 1{2n})\quad\text{ for all }0<\eps<\eps_{1n}.
$$ 

Let  the steady state $w_{AB}$ be defined by  \eqref{w_AB}.
Since $u_{\eps}$ is uniformly bounded in $Q$ there exists $M>0$ such that $u_\eps>-M$ in $Q$ for all $\eps\in (0,1]$.
It follows from  \eqref{w_A'} that for all $n\in \N$ there exists $A_n>1$ such that, independently of the choice of $B>0$,
$$
w_{A_nB}(y)<-M-1 \quad\text{for all }y\le -\tfrac 1{2n}.
$$
On the other hand there exist $y_n\in (\tfrac 1{2n},\tfrac 1n)$ and $B_n\in (0,1)$ such that 
$$
\tfrac 1n-y_n <\tfrac 12 \log \tfrac {1+B_n}{1-B_n}<\tfrac 12\left(\tfrac 12 (x_i+x_{i+1})-\tfrac 1n\right), \qquad w_{A_nB_n}(\tfrac 12 \log \tfrac {1+B_n}{1-B_n})<m_n
$$
(we recall that $\tfrac 12\log \tfrac {1+B_n}{1-B_n}$ is the point where $w_{A_nB_n}$ attains its maximal value). 
Finally we choose a point $z_n\in (\tfrac 12\log \tfrac {1+B_n}{1-B_n}, \log \tfrac {1+B_n}{1-B_n})$, which implies that 
$w_{A_nB_n}>0$ and $w'_{A_nB_n}<0$ at $z_n$. 

Let $w_{A_nB_n,\eps}$ be defined by Theorem \ref{thm: travelling waves approximation}.
By the choice of $A_n$, $B_n$  
and $z_n$ and by Theorem \ref{thm: travelling waves approximation}, there exists
$\eps_{2n}\le \eps_{1n}$ such that for all $\eps\in (0,\eps_{2n})$
$$
w_{A_nB_n,\eps }(y)<-M\text{ for }y\le -\tfrac 1{2n}, \quad 0<w_{A_nB_n,\eps }<m_n \text{ in } (0,z_n], \quad w'_{A_nB_n,\eps}(z_n)<0,
$$
and
$$
c_\eps T=\frac{(B_n-A_n)T}{2\log\varepsilon}<\frac 1n-y_n.
$$
We define $\underline u_{\eps,n}\in C([x_i,x_{i+1}]\times [0,T])$ by 
$$
\underline u_{\eps,n}(x,t)\!=\!
\begin{cases}
w_{A_nB_n,\eps }(x\!-\!x_i\!-\!y_n\!-\!c_\eps t) &\text{if }x_i\le x \le x_i\!+\!y_n\!+\!c_\eps t\!+\!z_n, \, 0\le t\le T\\
w_{A_nB_n,\eps }(x_{i+1}\!-\!y_n\!-\!c_\eps t\!-\!x)&\text{if }x_{i+1}\!-\!y_n\!-\!c_\eps t\!-\!z_n\le x \le x_{i+1}, \, 0\le t\le T\\
w_{A_nB_n,\eps }(z_n)&\text{otherwise.} 
\end{cases}
$$
By construction $\underline u_{\eps,n}$ is a (weak) subsolution of the parabolic equation for $u_\eps$ in $(x_i,x_{i+1})\times (0,T]$. In addition
$\underline u_{\eps,n}<-M\le u_\eps$ at $\{x_i,x_{i+1}\}\times (0,T)$ and $\underline u_{\eps,n}(x,0)\le u_{0\eps}(x)$ for $x\in (x_i,x_{i+1})$.
Hence $u_\eps\ge \underline u_{\eps,n}$ in $[x_i,x_{i+1}]\times (0,T]$. Since 
$$ 
x_i+y_n+c_\eps T<x_i+\tfrac 1n \quad\text{and}\quad  x_{i+1}-\tfrac 1n<x_{i+1}-y_n-c_\eps T
$$
and $\underline u_{\eps,n}(x,t)>0$ if $x_i+y_n+c_\eps T<x<x_{i+1}-y_n-c_\eps T$, this implies  that $u_\eps(x,t)>0$ if $x=x_i+\tfrac 1n$ or $x=x_{i+1}-\frac 1n$.

%%%%%%%%%%%%%%%%%%%%%%%%%%%%%%%%%%%%%%%%%%%%%%%%%%%%%
%%%%%%%%%%%%%%%%%%%%%%%%%%%%%%%%%%%%%%%%%%%%%%%%%%%%%
%%%%%%%%%%%%%%%%%%%%%%%%%%%%%%%%%%%%%%%%%%%%%%%%%%%%%
%%%%%%%%%%%%%%%%%%%%%%%%%%%%%%%%%%%%%%%%%%%%%%%%%%%%%
%%%%%%%%%%%%%%%%%%%%%%%%%%%%%%%%%%%%%%%%%%%%%%%%%%%%%
%%%%%%%%%%%%%%%%%%%%%%%%%%%%%%%%%%%%%%%%%%%%%%%%%%%%%
%%%%%%%%%%%%%%%%%%%%%%%%%%%%%%%%%%%%%%%%%%%%%%%%%%%%%
%%%%%%%%%%%%%%%%%%%%%%%%%%%%%%%%%%%%%%%%%%%%%%%%%%%%%
%%%%%%%%%%%%%%%%%%%%%%%%%%%%%%%%%%%%%%%%%%%%%%%%%%%%%
%%%%%%%%%%%%%%%%%%%%%%%%%%%%%%%%%%%%%%%%%%%%%%%%%%%%%
%%%%%%%%%%%%%%%%%%%%%%%%%%%%%%%%%%%%%%%%%%%%%%%%%%%%%
%%%%%%%%%%%%%%%%%%%%%%%%%%%%%%%%%%%%%%%%%%%%%%%%%%%%%
%%%%%%%%%%%%%%%%%%%%%%%%%%%%%%%%%%%%%%%%%%%%%%%%%%%%%

\section{Boundary regularity of the solution of the limit problem}\label{section regularity}

In this section we prove Theorem \ref{main result 2}, which exclusively concerns the unique solution $u$ of the limit problem \eqref{limit Dirichlet problem}.

\subsection{Proof of Theorem \ref{main result 2}$(i)$}
We show first that it is sufficient to prove the following lemma, 
which reminds the well-known estimate by Aronson and Benilan \cite{AB}  for nonnegative solutions
of the porous medium equation $u_t=\Delta (u^m)$ if $m>1$ (see also \cite{bertsch-discontinous}  in the case of nonnegative solutions of 
$u_t=u\Delta u -\gamma |\nabla u|^2$).

\begin{lem}
\label{lmm: bound on u_t}
Let $x_i,x_{i+1}\in \mathcal N_0$ and let $u_i$ be the restriction of $u$ to $\overline Q_i$, where we have set $Q_i=(x_i,x_{i+1})\times (0,T]$.  
Then 
$$
\text{$u_i>0$ in }Q_i\ \Rightarrow \ u_{it}\ge -\frac {u_i}t \quad\text{and}\quad u_{ixx}\ge -(1-u_i) -\frac 1t\quad\text{in }Q_i,
$$
$$
\text{$u_i<0$ in }Q_i\ \Rightarrow \ u_{it}\le -\frac {u_i}t \quad\text{and}\quad u_{ixx}\le -(1-u_i) -\frac 1t\quad\text{in }Q_i.
$$
\end{lem}

Consider for example the case that $u_i>0$. We recall that $u_{ix}$ is continuous in $Q_i$. 
By Lemma \ref{lmm: bound on u_t}, $u_{ixx}$ is bounded from below in $(x_i,x_{i+1})\times [t_0,T]$ for all 
$t_0\in (0,T]$. Since $u_i$ is bounded, this implies that the following one-sided limits are well defined: 
$$
u_{ix}(x_i^+,t)\in [-\infty,\infty),\qquad u_{ix}(x_{i+1}^-,t)\in (-\infty,+\infty].
$$ 
Since the function $x\mapsto u_i(x,t)$ attains a minimum (=0) at $x_i$ and $x_{i+1}$, we conclude that the limits $u_{ix}(x_i^+,t)$ and $u_{ix}(x_{i+1}^-,t)$ are 
bounded in $[t_0,T]$ for $t_0\in (0,T)$. If $u_i<0$ one argues similarly to arrive at the same conclusion. 
To complete the proof of part $(i)$ it remains to prove  Lemma \ref{lmm: bound on u_t}.

\medskip

\noindent {\it Proof of Lemma \ref{lmm: bound on u_t}.}
The inequalities for $u_{it}$ and $u_{ixx}$ in Lemma \ref{lmm: bound on u_t} are equivalent. Below we prove the inequality for $u_{it}$ in the case that $u_i>0$ in $Q_i$.

As we shall see in Section \ref{pf limit pb}, $u_i$ can be approximated from above by smooth solutions $u_{i,n}\ge \tfrac 1n$ 
of the uniformly parabolic problem 
\begin{equation}\label{approximate Dirichlet problem in Q_i bis}
\begin{cases}
u_t = u (u_{xx} + 1-u) &\text{in }Q_i\\
u(x_i,t)=u(x_{i+1},t)=\tfrac 1n &\text{for }t\in(0,T]\\
u(x,0)=u_0(x) +\tfrac 1n&\text{for }x\in ( x_i,x_{i+1}),
\end{cases}
\end{equation}
where $n\in \N$.
 We set
$$
p=\frac {(u_{i,n})_t}{u_{i,n}}\quad\text{in }Q_i.
$$
Then $(u_{i,n})_t=pu_{i,n}$ and $p=(u_{i,n})_{xx}+1-u_{i,n}$, whence
$$
\begin{aligned}
p_t&=(u_{i,n})_{txx}-(u_{i,n})_t=(pu_{i,n})_{xx}-pu_{i,n}=u_{i,n}p_{xx}+2(u_{i,n})_xp_x+((u_{i,n})_{xx}-u_{i,n})p\\
&=u_{i,n}p_{xx}+2(u_{i,n})_xp_x+(p-1)p\qquad\text{in }Q_i.
\end{aligned}
$$
Let $t_0\in (0,T)$ Since $p=0$ at $x_i$ and $x_{i+1}$ and since $\underline p=-(t-t_0)^{-1}$ is a subsolution of the equation in $(x_i,x_{i+1})\times (t_0,T)$
which tends to $-\infty$ as $t\to t_0^+$,
it follows easily from the Comparison Principle that $p\ge \underline p$, i.e.~$(t-t_0)(u_{i,n})_t\ge -u_{i,n}$, in $(x_i,x_{i+1})\times (t_0,T)$. 
Since $u_{i,n}\to u_i$ in $C^{2,1}_{\rm loc}(Q_i)$ as $n\to\infty$, this implies that  
$(t-t_0)u_{it}\ge -u_i$ in $(x_i,x_{i+1})\times (t_0,T)$. Since $t_0>0$ is arbitrary, we have proved that $tu_{it}\ge -u_i$ in $Q_i$.

\subsection{Proof of Theorem \ref{main result 2}$(ii)$}
Let $u_i$ and $Q_i$ be as above. We only consider the limit $u_i(x_i^+,t)$ in the case that $u_i>0$ in $Q_i$.
We define the difference quotient
$$
q(x,t)=\frac {u_i(x,t)-u_i(x_i,t)}{x-x_i}=\frac {u_{i}(x,t)}{x-x_i} \quad \text{for }t\in (t_0,T].
$$
Since, by part $(i)$, $u_{it}\ge -u_i/t$ in $Q_i$, it follows at once that also $q_t\ge -q/t$ in $Q_i$. Integration with respect to $t$ yields that
\begin{equation}\label{q estimate}
q(x,t)\ge \frac{t_0}{t}q(x,t_0) \quad \text{if }0< t_0<t\le T,
\end{equation}
whence also
$$
u_{ix}(x_i^+,t)\ge \frac{t_0}{t}u_{ix}(x_i^+,t_0)\quad \text{for }t\in (t_0,T].
$$
In particular $u_{ix}(x_i^+,t)>0$ if $u_{ix}(x_i^+,t_0)>0$ and $t>t_0$,
and we have proved  the existence of $\tau_i^+\in [0,\infty]$.

The function $t\mapsto u_{ix}(x_i^+,t)$ is continuous in $(0,\tau_i^+)$ since $u_{ix}(x_{i}^+,t)=0$ if $t\in (0,\tau_i^+)$. To prove the continuity 
in $(\tau_i^+,T]$ we shall show that
\begin{equation}\label{q cont}
q: (x_i,x_{i+1}) \times (\tau_i^+,T]\to \R \text{ can be extended by continuity to }  [x_i,x_{i+1}) \times (\tau_i^+,T].
\end{equation}

The function $q$ satisfies
\begin{equation*}
q_t = \frac{u_{it}}{x-x_i} = q u_{ixx} + q (1-u_i) =q ((x-x_i)q_{xx} +2q_x) +q(1-(x-x_i)q) \text{ in } (x_i,x_{i+1}) \times (\tau_i^+,T].
\end{equation*}
Setting $y= \sqrt{x-x_i}$ and $h(y,t)= q(x_i+y^2,t)$, we obtain that
\begin{equation*}
h_t = \frac{h}{4} \left(h_{yy} + 3 \frac{h_y}{y} \right) + h(1-y^2h) \text{ in } (0,\sqrt{x_{i+1}-x_i}) \times (\tau_i^+,T].
\end{equation*}
Let $B$ be the open ball in $\R^4$ centered at the origin with radius $\sqrt{x_{i+1}-x_i}$. Setting 
$y = |z|$ for $z\in B$ and $v(z,t) = h(|z|,t)$ for $(z,t)\in B\times (\tau_i^+,T]$, this means that
\begin{equation*}
v_t = \frac{v}{4} \Delta v + v(1-|z|^2 v) \text{ in } B \times (\tau_i^+,T].
\end{equation*}
It follows from \eqref{q estimate} and the definition of $\tau_i^+$ that, locally in $B \times (\tau_i^+,T]$, $v$ is bounded away from zero, whence 
locally in $B \times (\tau_i^+,T]$ the equation for $v$ 
is  uniformly parabolic. This implies the continuity of $v$ and we have proved \eqref{q cont}. This completes the proof of Theorem \ref{main result 2}.

\begin{remark}\label{continuity u_x}
It is not difficult to show that not only the functions $t\mapsto u_x(x_i^\pm,t)$ are continuous in $(0,T]\setminus\{\tau_i^\pm\}$, but also
the restriction of $u_x$ to $(x_i,x_{i+1})\times (0,T]$ can be extended with continuity to the set  
$[x_i,x_{i+1}]\times (0,T]\setminus \{(x_i,\tau_i^+),(x_{i+1},\tau_{i+1}^-)\}.$
\end{remark}

\begin{remark}\label{waiting time} 
The number $\tau^\pm_i$ reminds the concept of {\it waiting time} for the interfaces of the porous medium equation,
which is always a finite number \cite{vazquez-pme}. 
It is natural to ask whether also in our case $\tau_i^\pm$ is always finite. The answer is negative, as the following example shows.

Let $u_0$ be strictly increasing in $(a,b)$, let $u_0(x_1)=0$ and let $\log u_0\not\in L^1(x_1,b)$. 
Let $\psi\in C^1_{\rm c} ([x_1,b))$ be such that $\psi(x_1)=1$ and $\psi'\le 0$ in $(x_1,b)$ and set
$$
\chi_n(x)=\begin{cases} 
0 &\text{if }x_1\le x\le x_1+\tfrac 1n \\ 
n(x-x_i-\tfrac 1n )&\text{if }x_1+\tfrac 1n< x<x_1+\tfrac 2n \\
1 &\text{if }x_1+\tfrac 2n \le x\le b. \\ 
 \end{cases}
$$
Then $u_x>0$ in $(x_1,b)\times (0,T)$ and
$$
\begin{aligned}
&\int _{x_1}^b\log u(x,t) \psi(x)\chi_n(x) \dd x\\
&\qquad =\int _{x_1}^b\log u_0(x) \psi(x)\chi_n(x) \dd x+\iint _{(x_1,b)\times (0,t)}\left(-u_x(\psi'\chi_n+\psi\chi'_n)+(1-u)\psi\chi_n\right)\\
&\qquad \le\int _{x_i}^b\log u_0(x) \psi(x)\chi_n(x) \dd x+\iint _{(x_1,b)\times (0,t)}\left(-u_x\psi'\chi_n+(1-u)\psi\chi_n\right).\\
\end{aligned}$$
Letting $n\to\infty$ we obtain that 
$$
\int _{x_1}^b\!\log u(x,t) \psi(x) \dd x =-\infty \quad \text{for all }t\in (0,T],
$$
and since $T$ is arbitrary this implies that $\tau^+_1=\infty$.
\end{remark}

%%%%%%%%%%%%%%%%%%%%%%%%%%%%%%%%%%%%%%%%%%%%%%%%%%%%%
%%%%%%%%%%%%%%%%%%%%%%%%%%%%%%%%%%%%%%%%%%%%%%%%%%%%%
%%%%%%%%%%%%%%%%%%%%%%%%%%%%%%%%%%%%%%%%%%%%%%%%%%%%%
%%%%%%%%%%%%%%%%%%%%%%%%%%%%%%%%%%%%%%%%%%%%%%%%%%%%%
%%%%%%%%%%%%%%%%%%%%%%%%%%%%%%%%%%%%%%%%%%%%%%%%%%%%%
%%%%%%%%%%%%%%%%%%%%%%%%%%%%%%%%%%%%%%%%%%%%%%%%%%%%%
%%%%%%%%%%%%%%%%%%%%%%%%%%%%%%%%%%%%%%%%%%%%%%%%%%%%%
%%%%%%%%%%%%%%%%%%%%%%%%%%%%%%%%%%%%%%%%%%%%%%%%%%%%%
%%%%%%%%%%%%%%%%%%%%%%%%%%%%%%%%%%%%%%%%%%%%%%%%%%%%%
%%%%%%%%%%%%%%%%%%%%%%%%%%%%%%%%%%%%%%%%%%%%%%%%%%%%%
%%%%%%%%%%%%%%%%%%%%%%%%%%%%%%%%%%%%%%%%%%%%%%%%%%%%%
%%%%%%%%%%%%%%%%%%%%%%%%%%%%%%%%%%%%%%%%%%%%%%%%%%%%%
%%%%%%%%%%%%%%%%%%%%%%%%%%%%%%%%%%%%%%%%%%%%%%%%%%%%%

\section{The interface condition}\label{proof of main results}

This section is devoted to the proof of Theorem \ref{main result 3} concerning the asymptotic formula for the interface condition \eqref{conjecture}. 

The fact that $u_{\eps x}>0$ in $Q$ follows from the strong maximum principle applied to the equation for $u_{\eps x}$.
Hence the interface $x=\zeta_\eps(t)$ is well-defined by $u_\eps(\zeta_\eps(t),t)=0$ and, by the implicit function theorem, $\zeta_\eps\in C([0,T])\cap C^1((0,T])$.
By Theorem \ref{main result 1}, in particular by \eqref{interfaces do not move},
$$
\zeta_\eps\to x_1 \quad\text{in $C([0,T])$ as }\eps \to 0.
$$ 

For the limit function $u$ we have a similar positivity result of $u_x$, but we must exclude the point $x_1$ where the parabolic equation degenerates:
\begin{equation}\label{u_x positive}
u_{x}>0 \quad\text{in }  Q\setminus \{(x_1,t);\, t\in (0,T]\}.
\end{equation}

We begin with some preliminary calculations. Let $X_\eps:[-u_{1\eps},u_{1\eps}]\times [0,T]\to[a,b]$ be defined by 
$$
u_\eps(X_\eps(u,t),t)=u \qquad (\text{and so }\zeta_\eps(t)=X_\eps(0,t)).
$$   
Differentiating $x = X_\varepsilon(u,t)$ with respect to $x$ and $t$ we find that  
$$
X_{\varepsilon u }(u,t) = \frac{1}{u_{\varepsilon x} (X_\varepsilon(u,t),t)}, \qquad 
X_{\varepsilon t} (u,t) = -\frac{u_{\varepsilon t} (X_\varepsilon(u,t),t)}{u_{\varepsilon x} (X_\varepsilon(u,t),t)},
$$
and it follows from the equation for $u_\eps$ that $X_\eps$ satisfies the parabolic equation
\begin{equation}\label{eq: equation with eps after change of variable X}
X_t = -(\varepsilon+\Phi^2_\varepsilon(u))\left(\frac{1}{X_u} \right)_u -  \Phi_\varepsilon(u)(1-\Phi^2_\varepsilon(u))\sqrt{\varepsilon+\Phi^2_\varepsilon(u)} X_u.
\end{equation}

We define $X_0:[-1,1]\times [0,T]\to[a,b]$ by 
$$
u(X_0(u,t),t)=u \qquad (\text{and so }X_0(0,t)=x_1),
$$   
and one shows in a similar way that $X_0$ satisfies
\begin{equation}\label{eq: equation X_0}
X_t = -|u|\left(\frac{1}{X_u} \right)_u -  u(1-|u|) X_u \quad\text{in } ([-1,0)\cup(0,1])\times (0,T].
\end{equation}
It easily follows from Theorem \ref{main result 1} and \eqref{u_x positive} that 
$$
X_\eps\to X_0 \quad\text{in } C([-1,1]\times[0,T])\cap C^{2,1}_{\rm loc}(([-1,0)\cup(0,1])\times (0,T]).
$$

Since, by Theorem \ref{main result 1}, $u_{\eps x}\to u_x$ locally in $([a,b]\setminus\{x_1\})\times (0,T]$, for all $\eps\in (0,1]$ there exists $0<\delta_\eps\to 0$ such that 
$$
\frac 1{X_{\eps u}(\pm\delta_\eps,\cdot)}-\frac 1{X_{0 u}(\pm\delta_\eps,\cdot)}\to 0 \quad\text{in 
$C_{\rm loc}((0,T])$ as }\eps\to 0.
$$
Hence it follows from the  continuity properties of $u_x$ (see Remark \ref{continuity u_x}) that
\begin{equation}\label{choice delta}
\frac 1{X_{\eps u}(\pm\delta_\eps,\cdot)}\to u_x(x_1^\pm,\cdot) \quad\text{pointwise in $(0,T]$ and in 
$C_{\rm loc}((0,T]\setminus \{\tau^\pm_1\})$ as }\eps\to 0.
\end{equation}

We fix $t\in (0,T]$.
Without loss of generality we may assume that $\eps/\delta_\eps\to 0$ as $\eps\to 0$.
The idea of the proof is to integrate the equation for $X_\varepsilon$ with respect to $u$ in a neighbourhood of $u=0$, but to do so we change variable and set 
\begin{equation*}
p = A_\varepsilon(u) := \int_0^u \frac{1}{\varepsilon+\Phi_\varepsilon^2(s)}\dd{s}.
\end{equation*}
Then the second order term in \eqref{eq: equation with eps after change of variable X} becomes a partial derivative with respect to $p$ and
integrating  \eqref{eq: equation with eps after change of variable X}  with respect to $p$ from $-A_\varepsilon(\delta_\eps)$ to $A_\varepsilon(\delta_\eps)$ we obtain that
\begin{equation}\label{I,II,III}
\begin{aligned}
&\int_{-A_\varepsilon(\delta_\eps)}^{A_\varepsilon(\delta_\eps)} X_{\varepsilon t} (A_\varepsilon^{-1} (p),t) \, \dd{p} 
+\frac 1{X_{\eps u}(\delta_\eps,t)}-\frac 1{X_{\eps u}(-\delta_\eps,t)}\\ 
&\qquad \qquad =-B_{\varepsilon,\delta_\eps}(t):=- \int_{-\delta_\eps}^{\delta_\eps} \frac{\Phi_\varepsilon(u)(1-\Phi_\varepsilon^2(u))}{\sqrt{\varepsilon+\Phi_\varepsilon^2(u)}} X_{\varepsilon u}(u,t) \, \dd{u}.
\end{aligned}
\end{equation}
Observe that, as $\eps\to 0$, 
$$
|B_{\varepsilon,\delta_\eps}(t)| \le  \int_{-\delta_\eps}^{\delta_\eps}  X_{\varepsilon u}(u,t) \, \dd{u}=X_\eps(\delta_\eps,t)-X_\eps(-\delta_\eps,t)\to 0 
$$
uniformly with respect to $t$; here we have used the (uniform) continuity of the map $[0,1]\times [-\tfrac 12 ,\tfrac12]\times [0,T]\ni(\eps,u,t)\mapsto X_\eps(u,t)$. 
In view of \eqref{choice delta} and \eqref{I,II,III} this implies that
\begin{equation}\label{FBC 1}
-\int_{-A_\varepsilon(\delta_\eps)}^{A_\varepsilon(\delta_\eps)} X_{\varepsilon t} (A_\varepsilon^{-1} (p),t) \, \dd{p}\to  u_x(x_1^+,t)-u_x(x_1^-,t)\quad\text{as }\eps\to0.
\end{equation}

We claim that 
\begin{equation}\label{A eps = -log eps}
A_\eps(\delta_\eps) = -\log \eps (1+o(1))\quad\text{if $\delta_\eps \to 0$ and $\eps = o(\delta_\eps)$ as $\eps \to 0$.}
\end{equation}
By the definition of $\Phi_\varepsilon$,
\begin{equation*}
u = 2\int_0^{\Phi_\varepsilon(u)} \sqrt{\varepsilon+s^2} \, \dd{s}\ \Rightarrow\ 1 = 2\sqrt{\varepsilon+\Phi_\varepsilon^2(u)} \Phi_\varepsilon'(u),
\end{equation*}
whence, setting $\phi = \Phi_\varepsilon(s)$, 
\begin{equation*}
A_\varepsilon(\delta_\eps) \!=\! \int_0^{\delta_\eps}\!\! \frac{1}{\varepsilon \!+\! \Phi_\varepsilon^2(s)}\dd{s} \!=\!2\!\int_0^{\Phi_\varepsilon(\delta_\eps)} 
\!\!\! \frac{1}{\sqrt{\varepsilon\!+\!\phi^2}}\dd{\phi}
\!=\! 2\log \left(\frac{\Phi_\varepsilon(\delta_\eps)}{\sqrt{\varepsilon}}\!+\!\sqrt{1\!+\!\left( \frac{\Phi_\varepsilon(\delta_\eps)}{\sqrt{\varepsilon}}\right)^2}  \right).
\end{equation*} 
Recall that $\Phi_\eps = U_\eps^{-1}$ with $U_\eps$ defined in \eqref{eq: change of variable}:
$$
U_\eps(\phi) = 2\int_0^{\phi} \sqrt{\varepsilon+s^2} \, \dd{s} = \phi \sqrt{\varepsilon+\phi^2}
+\varepsilon \log\left(\frac \phi{\sqrt\varepsilon}+\sqrt{1+\frac{\phi^2}\varepsilon}\right),
$$
and observe that $U_\eps(\sqrt\eps)= (\sqrt 2+\log(1+\sqrt 2))\eps$.
Setting $\phi_\eps= \Phi_\eps(\delta_\eps)$ and $\xi_\eps = \frac{\phi_\eps}{\sqrt{\eps}}$,
$$
\frac{\delta_\eps}\eps=\frac{U_\eps(\phi_\eps)}\eps=\frac{U_\eps(\xi_\eps\sqrt \eps)}\eps
= \xi_\eps \sqrt{1+\xi_\eps^2}+\log\left(\xi_\eps+\sqrt{1+\xi_\eps^2}\right),
$$
and since $\delta_\eps/\eps\to \infty$ as $\eps\to0$ this implies that 
$\xi_\eps\to\infty$ as $\eps\to0$. Hence
$$
\frac{\Phi_\eps(\delta_\eps)}{\sqrt{\delta_\eps}} = \frac{\phi_\eps}{\sqrt{U_\eps(\phi_\eps)}} = 
\frac{\xi_\eps}{\sqrt{\xi_\eps \sqrt{1+\xi_\eps^2} + \log\left( \xi_\eps + \sqrt{1+ \xi_\eps^2} \right) }} \to 1 \text{ as } \eps \to 0
$$
and we obtain \eqref{A eps = -log eps}: as $\eps\to 0$
$$
A_\eps(\delta_\eps) = 2 \log \left( 2\frac{\sqrt{\delta_\eps}}{\sqrt{\varepsilon}}(1+o(1))\right) = 2\log(2\sqrt{\delta_\eps})-\log\varepsilon + o(1) %\\ &
= -\log \varepsilon (1+o(1)). 
$$

It follows from \eqref{FBC 1} and \eqref{A eps = -log eps} that
\begin{align*}
\left( \int_{-\delta_\eps}^{\delta_\eps} \frac{\dd{u}}{\varepsilon+\Phi_\varepsilon^2(u)} \right)^{-1}  \int_{-\delta_\eps}^{\delta_\eps}  \frac{2\log \eps X_{\eps t}(u,t)}{\varepsilon+\Phi_\varepsilon^2(u)} \dd{u}
 &= \frac {\log \eps}{A_\eps{(\delta_\eps)}}  \int_{-A_\varepsilon(\delta_\eps)}^{A_\varepsilon(\delta_\eps)} X_{\varepsilon t} (A_\varepsilon^{-1} (p),t) \, \dd{p}
 \\&
 \to  u_x(x_1^+,t)-u_x(x_1^-,t)\quad\text{as }\eps\to0.
\end{align*}
This completes the proof of Theorem \ref{main result 3}.

\begin{remark}\label{possible weak formulations}
The weak formulation of \eqref{conjecture}, which is given by Theorem \ref{main result 3} and proved in this Section, is not the only possible weak version of the
asymptotic formula for the velocity  of the interface, $\zeta'_\eps(t)$.  For example, if
$$
\tau_1=\max\{\tau_1^+,\tau_1^-\}<\infty,
$$
one can show that
\begin{equation}\label{weak FBC}
2\log \eps\,  \zeta_\eps' \text{ converges weakly in $L^2_{\rm loc}(\tau_1,T)$ to }u_x(x_1^+,\cdot)-u_x(x_1^-,\cdot)\text{ as }\eps\to0.
\end{equation}

Indeed, in view of Theorem \ref{main result 3} we obtain \eqref{weak FBC} if we  prove that, given  $\varphi \in  C^\infty_c((\tau,T))$, 
\begin{align*}
&\lim_{\eps\to 0} 2\log \eps \left( \int_{-\delta_\eps}^{\delta_\eps} \frac{\dd{u}}{\varepsilon+\Phi_\varepsilon^2(u)} \right)^{-1} \int_0^T \int_{-\delta_\eps}^{\delta_\eps}  \frac{X_{\eps t}(u,t) \varphi(t)}{\varepsilon+\Phi_\varepsilon^2(u)} \dd{u} \dd{t} 
\\&\qquad
= \lim_{\eps\to 0} 2 \log \eps \int_0^T X_{\eps t}(0,t) \varphi(t) \dd{t}.
\end{align*}
Due to \eqref{A eps = -log eps}, this is equivalent to proving that  
$$
\lim_{\eps \to 0} \int_{-\delta_\eps}^{\delta_\eps} \frac{1}{\eps+\Phi_\eps^2(u)} \left( \int_0^T (X_{\eps t}(u,t) - X_{\eps t}(0,t)) \varphi(t) \dd{t}  \right) \dd{u} = 0.
$$
By the definition of $\tau_1^\pm$, there exists $C>0$ which does not depend on $\eps$ such that $|X_{\eps u}(x,t) |< C$ for $x\in (a,b)$ and $t\in \supp \varphi$.
Hence 
$$
\begin{aligned}
&\left| \int_{-\delta_\eps}^{\delta_\eps} \frac{1}{\eps+\Phi_\eps^2(u)} \left( \int_0^T (X_{\eps t}(u,t) - X_{\eps t}(0,t)) \varphi(t) \dd{t}  \right) \dd{u}\right| \\
&\qquad =\left|  \int_{-\delta_\eps}^{\delta_\eps} \frac{1}{\eps+\Phi_\eps^2(u)} \left( \int_0^T \left(\int_0^uX_{\eps u }(s,t) \dd s\right) \varphi'(t) \dd{t}  \right) \dd{u}\right|\\
&\qquad \le CT\|\varphi'\|_\infty\int_{-\delta_\eps}^{\delta_\eps} \frac{|u|}{\eps+\Phi_\eps^2(u)}     \dd{u}\to 0 \quad\text{as }\eps \to 0,\\
\end{aligned}
$$
since 
$$
\frac{|u|}{\eps + \Phi_\eps^2(u)} \leq \frac{|u|}{\Phi_\eps^2(u)} \to 1 \quad\text{as }\eps \to 0.
$$

\end{remark}

%%%%%%%%%%%%%%%%%%%%%%%%%%%%%%%%%%%%%%%%%%%
%%%%%%%%%%%%%%%%%%%%%%%%%%%%%%%%%%%%%%%%%%%
%%%%%%%%%%%%%%%%%%%%%%%%%%%%%%%%%%%%%%%%%%%
%%%%%%%%%%%%%%%%%%%%%%%%%%%%%%%%%%%%%%%%%%%
%%%%%%%%%%%%%%%%%%%%%%%%%%%%%%%%%%%%%%%%%%%
%%%%%%%%%%%%%%%%%%%%%%%%%%%%%%%%%%%%%%%%%%%
%%%%%%%%%%%%%%%%%%%%%%%%%%%%%%%%%%%%%%%%%%%
%%%%%%%%%%%%%%%%%%%%%%%%%%%%%%%%%%%%%%%%%%%
%%%%%%%%%%%%%%%%%%%%%%%%%%%%%%%%%%%%%%%%%%%
%%%%%%%%%%%%%%%%%%%%%%%%%%%%%%%%%%%%%%%%%%%
%%%%%%%%%%%%%%%%%%%%%%%%%%%%%%%%%%%%%%%%%%%
%%%%%%%%%%%%%%%%%%%%%%%%%%%%%%%%%%%%%%%%%%%
%%%%%%%%%%%%%%%%%%%%%%%%%%%%%%%%%%%%%%%%%%%
%%%%%%%%%%%%%%%%%%%%%%%%%%%%%%%%%%%%%%%%%%%

\section{The limit problem: proof of Theorem \ref{well-posedness limit problem}}\label{pf limit pb}

Consider the Dirichlet problem in $( x_i,x_{i+1})\times (0,T]$:
\begin{equation}\label{Dirichlet problem in Q_i}
\begin{cases}
u_t = |u| u_{xx} + u(1-|u|) &\text{in }Q_i:=( x_i,x_{i+1})\times (0,T]\\
u(x_i,t)=u(x_{i+1},t)=0 &\text{for }t\in(0,T]\\
u(x,0)=u_0(x) &\text{for }x\in ( x_i,x_{i+1}).
\end{cases}
\end{equation}
We assume that $u_0>0$ in $(x_i,x_{i+1})$ (the case that $u_0<0$ is completely similar) and approximate problem \eqref{Dirichlet problem in Q_i} by 
\begin{equation}\label{approximate Dirichlet problem in Q_i}
\begin{cases}
u_t = u (u_{xx} + 1-u) &\text{in }Q_i\\
u(x_i,t)=u(x_{i+1},t)=\tfrac 1n &\text{for }t\in(0,T]\\
u(x,0)=u_0(x) +\tfrac 1n&\text{for }x\in ( x_i,x_{i+1}),
\end{cases}
\end{equation}
where $n\in \N$. We collect various basic results on problems \eqref{Dirichlet problem in Q_i} and \eqref{approximate Dirichlet problem in Q_i} in:

\begin{lem}\label{results limit problem}
Let $u_0$ satisfy \eqref{hyp u_0} and let $u_0>0$ in $(x_i,x_{i+1})$. Then problem \eqref{approximate Dirichlet problem in Q_i} has, for all $n\in\N$,  
a unique solution $u_{i,n}\in C(\overline Q_i)\cap C^{2,1}([x_i,x_{i+1}]\times (0,T])$, $u_{i,n}$ is pointwise decreasing with respect to $n$, and its pointwise 
limit $u_i$ is the unique solution in $C(\overline Q_i)\cap C^{2,1}(Q_i)$ of problem \eqref{Dirichlet problem in Q_i} which is positive in $Q_i$.
\end{lem}

\begin{proof}
By the Comparison Principle for uniformly parabolic equations, 
each smooth solution $u$ of \eqref{approximate Dirichlet problem in Q_i} satisfies $\tfrac 1n\le u\le \max\{1, \|u_0\|_\infty+\tfrac 1n\}$, whence,
by standard theory of quasilinear parabolic equations, problem \eqref{approximate Dirichlet problem in Q_i} has a unique solution 
$u_{i,n}$ in $ C(\overline Q_i)\cap C^{2,1}([x_i,x_{i+1}]\times (0,T])$. In addition $u_{i,n}$ is pointwise decreasing with respect to $n$ and there exists
$$
u_i(x,t)=\lim_{n\to \infty} u_{i,n}(x,t)\ge 0 \quad\text{for }(x,t)\in \overline Q_i.
$$

It follows easily from the construction of explicit positive subsolutions (not depending on $t$ and $n$) in subintervals of $(x_i,x_{i+1})$ 
that $u_{i,n}(x,t)\ge g_i(x)>0$ in $Q_i$ for some continuous function $g$, whence also $u_i(x,t)\ge g(x)>0$. 
Hence, by standard a priori local Schauder type bounds (\cite {LSU}) 
for solutions  of quasilinear parabolic equations, 
\begin{equation}\label{schauder}
u_{i,n}\to u_i\quad\text{in }C_{\text{loc}}^{2,1}(Q_i),
\end{equation}
and $u_i$ satisfies the PDE pointwise in $Q_i$. 
Since $u_{i,n}=\tfrac 1n$ at $x_i$ and $x_{i+1}$ and since $u_{i,n}$ is decreasing in $n$, $u_i$ vanishes and is continuous on 
the lateral boundary $\{x_i,x_{i+1}\}\times [0,T]$ of $Q_i$.  
In addition, by local H\"older estimates (\cite {LSU}) in $(x_i,x_{i+1})\times [0,T]$  and the Lipschitz continuity of $u_0$, $u_i$ is continuous in $(x_i,x_{i+1})\times\{0\}$.
Hence $u_i\in C(\overline Q_i)\cap C^{2,1}(Q_i)$ is a solution of problem \eqref{Dirichlet problem in Q_i} which is strictly positive in $Q_i$.

It remains to prove the uniqueness claim of Lemma \ref{results limit problem}.
Let  $v$ be another classical solution such that $v>0$ in $Q_i$. First we show that 
\begin{equation}\label{maximality}
0<v\le u_i\quad \text{in }Q_i.
\end{equation}
Let $n\in \N$ and let $\delta_n>0$ be so small that $v(x_i+\delta_n,t)<\tfrac 1n$
and $v(x_{i+1}-\delta_n,t)<\tfrac 1n$ for $t\in [0,T]$. Then it follows from the Comparison Principle for uniformly parabolic equations that $v<u_{i,n}$ in 
$(x_i+\delta_n,x_{i+1}-\delta_n)\times [0,T]$. Since $n$ is arbitrary and $\delta_n\to 0$ as $n\to \infty$, this implies \eqref{maximality}.

To show that $v=u_i$ in $Q_i$ we observe that, by the equations for $u_i$ and $v$, 
\begin{equation}\label{log}
(\log u_i - \log v)_t = (u_i-v)_{xx} - (u_i-v) \text{ in } Q_i.
\end{equation}
We would like to test this equation with the first eigenfunction $\varphi_0$ of the laplacian in $(x_i,x_{i+1})$
with homogeneous Dirichlet data such that $\max  \varphi_0=1$. But since \eqref{log} is singular at the lateral boundaries 
we slightly shrink the interval: let $\delta>0$ be sufficiently small and let $\varphi_\delta$ be the first eigenfunction of the laplacian in $(x_i+\delta,x_{i+1}-\delta)$ 
with homogeneous Dirichlet data such that $\max  \varphi_\delta=1$. We denote the first eigenvalue by $\lambda_\delta$:
$$
\lambda_\delta =-\frac{\pi^2}{(x_{i+1}-x_i-2\delta)^2}\to\lambda_0:=-\frac{\pi^2}{(x_{i+1}-x_i)^2}\quad \text{as }\delta \downarrow0. 
$$
In addition $\varphi_\delta \uparrow \varphi$ locally in $(x_i,x_{i+1})$ as $\delta\downarrow 0$.

Let $t \in (0,T]$ and $\delta>0$. Multiplying \eqref{log} by $\varphi_\delta$ and integrating, we obtain
\begin{equation*}
\begin{aligned}
&\int_{x_i+\delta}^{x_{i+1}-\delta}(\log u_i(t) - \log v(t)) \varphi_\delta \dd x+
(1-\lambda_\delta) \int_0^t \int_{x_i+\delta}^{x_{i+1}-\delta} (u_i-v) \varphi_\delta \, \dd{x}  \dd{t} \\
&\qquad  =\sqrt{|\lambda_\delta|} \left(\int_{0}^{t}(u_i(x_{i+1}-\delta,t)-v(x_{i+1}-\delta,t))\dd t+  \int_{0}^{t}(u_i(x_i+\delta,t)-v(x_i+\delta,t))\dd t\right).
\end{aligned}
\end{equation*}
Since the integrals on the right-hand side vanish as $\delta\to 0$, it follows from \eqref{maximality} and the monotone convergence theorem that
$$
0\le \int_{x_i}^{x_{i+1}}(\log u_i(t) - \log v(t)) \varphi_0 \dd x=
-(1-\lambda_0) \int_0^t \int_{x_i}^{x_{i+1}} (u_i-v) \varphi_0 \, \dd{x}  \dd{t}\le 0.
$$
Hence, by \eqref{maximality}, the continuity of $u_i$ and $v$ and the arbitrariness of $t$,  
we conclude that $u_i=v$ in $Q_i$.
\end{proof}

Lemma \ref{results limit problem} concerns the sets $Q_i$ for $i=1\dots, k-1$,
but a similar result can be easily proved  in $(a,x_1)\times (0,T]$ and $(x_k,b)\times (0,T]$. 
At this point the solutions in the (k+1) single intervals can be ``merged together'' to define a function $u$ in all of $Q$, 
i.e.~the restriction of $u$ to one of the sets $Q_i$ coincides with the smooth solution of problem \eqref{Dirichlet problem in Q_i} in $Q_i$. 
By construction, $u$ is a classical solution of problem \eqref{limit Dirichlet problem}.

Viceversa, the restriction of any classical solution of problem \eqref{limit Dirichlet problem} to one of the sets $Q_i$ is a smooth solution of 
problem \eqref{Dirichlet problem in Q_i} in $Q_i$, so the uniqueness statement in Lemma \ref{results limit problem} implies the uniqueness of the classical solution 
of problem \eqref{Dirichlet problem in Q_i}. 

To complete the proof of Theorem \ref{well-posedness limit problem} we must 
show that $u$ is a weak solution of problem \eqref{limit Dirichlet problem}.
This is an immediate consequence of the integral equality \eqref{weak u_i} in the following result.

\begin{lem}
\label{lmm: L^2 bound on gradient}
Let $u_{i,n}$ and $u_i$ be  defined by Lemma \ref{results limit problem} and let $\alpha > -1$.
Then there exists a positive constant $K$ which does not depend on $n$ such that 
\begin{equation}\label{gradient}
\tfrac{4(\alpha+1)}{(\alpha+2)^2}\iint_{Q_i} \left(u_{i,n}^{\frac{\alpha+2}{2}}\right)_x^2  \dd{x} \dd{t} 
+\!\int_0^T n^{-(\alpha+1)} \left(|(u_{i,n})_x(x_{i+1},t)| + |(u_{i,n})_x(x_i,t)|\right) \dd{t}\leq K.
\end{equation} 
In addition $u_i^{\frac {\alpha+2}2}\in L^2(0,T; H^1((x_i,x_{i+1})))$,
\begin{equation}\label{weak conv grad}
\frac{\partial}{\partial x}\left(u_{i,n}^{\frac {\alpha+2}2}\right)\rightharpoonup \frac{\partial}{\partial x}\left(u_i^{\frac {\alpha+2}2}\right)\quad\text{in $L^2(Q_i)$ as $n\to\infty$},
\end{equation}
and, for all $\psi \in C^{1,1}_c([x_i,x_{i+1}]\times [0,T))$,
\begin{equation}\label{weak u_i}
\int_{x_i}^{x_{i+1}} \!\!u_0(x) \psi(x,0)\dd x + \!\iint_{Q_i}\left(u_i \psi_t- u_iu_{ix}  \psi_x  -  u_{ix}^2 \psi +u_i(1-u_i)\psi\right) \dd x\dd t = 0 .
\end{equation} 
\end{lem}

\begin{proof}
Let $\tau\in (0,T)$. Integration by parts over $Q_i^\tau:=(x_i,x_{i+1})\times (\tau,T)$ yields 
$$
\begin{aligned}
&\tfrac{4(\alpha+1)}{(\alpha+2)^2}\iint_{Q_i^\tau} \left(u_{i,n}^{\frac{\alpha+2}{2}}\right)_x^2  \dd x \dd t 
-\!\int_\tau^T n^{-(\alpha+1)} \left((u_{i,n})_x(x_{i+1},t) -  (u_{i,n})_x(x_i,t)\right) \dd t\\
&\qquad 
= \frac{1}{\alpha+1} \int_{x_i}^{x_{i+1}} \left(u_{i,n}^{\alpha+1}(x,\tau) - u_{i,n}^{\alpha+1}(x,T)\right) \dd{x} 
+ \iint_{Q_{\tau,T}} u_{i,n}^{\alpha+1} (1-u_{i,n})  \dd{x} \dd{t}.
\end{aligned}
$$
Since $u_{i,n}=\tfrac 1n$ at $x_i$ and $x_{i+1}$ and $u_{i,n}\ge \tfrac 1n$ in $Q_i$, we have that
$(u_{i,n})_x(x_{i+1},t)\le 0$ and $(u_{i,n})_x(x_{i},t)\ge 0$ for $t\in (0,T]$. Letting $\tau\to 0$, the existence of the constant $K$ follows 
from the uniform boundedness of $u_{i,n}$ in $Q_i$.

It remains to prove \eqref{weak u_i}. Let $\psi \in C^{1,1}_c([x_i,x_{i+1}]\times [0,T))$ and let $\tau\in (0,T)$. 
Integration by parts over $Q_i^\tau$ yields
\begin{equation}\label{weak u_i bis}
\begin{aligned}
&\int_{x_i}^{x_{i+1}} \!\!u_{i,n}(x,\tau) \psi(x,\tau)\dd x  +\tfrac 1n\int_\tau^T  (u_{i,n})_x\psi \big|^{(x_{i+1},t)}_{(x_i,t)} \dd t\\
&\qquad 
+ \!\iint_{Q_i^\tau}\!\!\left(u_{i,n} \psi_{t}\!-\! \tfrac 12 (u^2_{i,n})_x \psi_{x}\!-\!(u_{i,n})_x^2 \psi \!+\!u_{i,n}(1\!-\!u_{i,n})\psi\right) \dd x\dd t = 0. 
\end{aligned}
\end{equation}

We first pass to the limit $n\to\infty$ inside the integrals, for fixed $\tau$.
By the Dominated Convergence Theorem, the only nontrivial terms are those containing $(u_{i,n})_x$. 
In the term with $(u_{i,n}^2)_x$ it is enough to use \eqref{weak conv grad} with $\alpha=2$ to pass to the limit.
In addition, we obtain from \eqref{gradient} with $\alpha\in (-1,0)$ that the integral 
$$
\tfrac 1n\left|\int_\tau^T  (u_{i,n})_x\psi \big|^{(x_{i+1},t)}_{(x_i,t)} \dd t\right|\le C n^\alpha \int_\tau^T n^{-(\alpha+1)} \left(|(u_{i,n})_x(x_{i+1},t)| + |(u_{i,n})_x(x_i,t)|\right)
$$
vanishes as $n\to \infty$. It remains to consider the term containing $(u_{i,n})_x^2$. By \eqref{schauder}, for each $\delta>0$
$$
\iint_{(x_i+\delta,x_{i+1}-\delta)\times (\tau,T)} (u_{i,n})_x^2 \psi \, \dd x\dd t \to\iint_{(x_i+\delta,x_{i+1}-\delta)\times (\tau,T)} (u_{i})_x^2 \psi \, \dd x\dd t
\quad\text{as }n\to\infty.
$$
On the other hand, by the pointwise monotonicity of $u_{i,n}$ with respect to $n$ and the lateral boundary condition $u_{i,n}=\tfrac 1n$, for all $\eta>0$ 
there exist $\delta_\eta>0$ and $n_\eta\in \N$ such that $u_{i,n}(x,t)<\eta$ if $n>n_\eta$ and if $x<x_i +\delta_\eta$ or  $x>x_{i+1} -\delta_\eta$.
Hence, by \eqref{gradient} with $\alpha=-1/2$,
$$
\begin{aligned}
&\left|\,\iint_{(x_i,x_i+\delta_\eta)\times (\tau,T)} (u_{i,n})_x^2 \psi \, \dd x\dd t\, \right|\\
&\qquad \le \left(\sup_{(x_i,x_i+\delta_\eta)\times (\tau,T)} \sqrt {u_{i,n}}\right)
\iint_{(x_i,x_i+\delta_\eta)\times (\tau,T)} u_{i,n}^{-1/2}(u_{i,n})_x^2 |\psi| \, \dd x\dd t\\
&\qquad \le C_1\sup_{(x_i,x_i+\delta_\eta)\times (\tau,T)} \sqrt {u_{i,n}}\le  C_1\sqrt \eta \quad \text{for all } n>n_\eta.
\end{aligned}
$$
Since $\eta$ is arbitrary this implies that we can let $n\to\infty$ in  \eqref{weak u_i bis}:
$$
\int_{x_i}^{x_{i+1}} \!\!u_{i}(x,\tau) \psi(x,\tau)\dd x  \\
+ \!\iint_{Q_i^\tau}\!\!\left(u_{i} \psi_{t}\!-\! u_i u_{ix} \psi_{x}\!-\!u_{ix}^2 \psi \!+\!u_{i}(1\!-\!u_{i})\psi\right) \dd x\dd t = 0. 
$$
Since $u_{ix}\in L^2(Q_i)$ we can let $\tau\to 0$ and obtain \eqref{weak u_i}.
\end{proof}

\bigskip
\noindent
{\bf Acknowledgement.}
MB acknowledges the MIUR Excellence Department Project awarded to the Department of Mathematics, University of Rome Tor Vergata, CUP E83C18000100006.

\end{document}